\documentclass[12pt]{article}

\usepackage{amsmath}
\usepackage{amssymb}

\numberwithin{equation}{section}

\newtheorem{thm}{Theorem}[section]

\newtheorem{cor}[thm]{Corollary}
\newtheorem{dfn}[thm]{Definition}

\newtheorem{lem}[thm]{Lemma}

\newcommand{\eqa}{\begin{eqnarray}}
\newcommand{\eeqa}{\end{eqnarray}}
\newcommand{\beq}{\begin{equation}}
\newcommand{\eeq}{\end{equation}}
\newcommand{\nn}{\nonumber}
\newcommand{\p}{\partial}
\newcommand{\ad}{{\mbox{ad}}}
\newcommand{\vp}{\varphi}
\newcommand{\ve}{\epsilon}
\newcommand{\dl}{\delta}

\newcommand{\lm}{\lambda}
\newcommand{\hx}{\hat X}
\newcommand{\hy}{\hat Y}

\newcommand{\pal}{\partial}

\newcommand{\al}{\alpha}
\newcommand{\bt}{\beta}

\newcommand{\tg}{\tilde{g}}
\newcommand{\tGamma}{\tilde{\Gamma}}
\newcommand{\tnabla}{\tilde{\nabla}}

\newcommand{\pf}{\noindent{\it Proof \ }}
\newcommand{\epf}{$\quad$\hfill
\raisebox{0.11truecm}{\fbox{}}\par\vskip0.4truecm}
\newenvironment{prf}{\noindent {\it Proof} \ }{\hfill $\Box$}

\begin{document}

\title {On Hamiltonian perturbations of hyperbolic systems of conservation laws}

\author{{Boris Dubrovin${}^*$ \ \ Si-Qi Liu${}^{**}$ \ \ Youjin Zhang${}^{**}$}\\
{${}^*$\small SISSA, Via Beirut 2-4, 34014 Trieste, Italy, and Steklov
Math. Institute, Moscow}\\
 {${}^{**}$\small Department of Mathematical Sciences, Tsinghua
University}\\
{\small Beijing 100084, P.R.China}}
\date{}
\maketitle

\begin{abstract}
We study the general structure of formal
perturbative solutions to the Hamiltonian
perturbations of spatially one-dimensional systems of hyperbolic PDEs
${\bf v}_t +[\phi({\bf v})]_x=0$. Under certain genericity assumptions
it is proved that any {\it bihamiltonian} perturbation can be eliminated in all
orders of the perturbative expansion by a change of coordinates on the infinite jet
space depending {\it rationally} on the derivatives. The main tools is in constructing
of the so-called quasi-Miura transformation of jet coordinates eliminating
an arbitrary
deformation of a semisimple bihamiltonian
structure of hydrodynamic type (the quasitriviality theorem). We also describe,
following \cite{LZ1}, the invariants of such bihamiltonian structures with
respect to the group of Miura-type transformations depending {\it polynomially}
on the derivatives.
\end{abstract}

\tableofcontents

\section{Introduction}

Systems of evolutionary PDEs arising in many physical
applications can be written in the form
\beq\label{pert1}
w^i_t +V^i_j(w) w^j_x + \mbox{perturbation}=0, \quad i=1, \dots, n
\eeq
where the perturbation may depend on higher derivatives.
The dependent variables of the system
$$
w=\left(w^1(x,t), \dots, w^n(x,t)\right)
$$
are functions of one spatial variable $x$ and the time $t$, summation over
repeated indices will be assumed.
$V^i_j(w)$ is a matrix of functions having {\it real distinct} eigenvalues.
Therefore
the system (\ref{pert1}) can be considered as
a perturbation of the hyperbolic system of first order quasilinear PDEs
\beq\label{pert2}
v^i_t + V^i_j(v) v^j_x =0, \quad i=1, \dots, n
\eeq
(it will be convenient to denote differently the dependent variables of the
unperturbed system (\ref{pert2}) and the perturbed one (\ref{pert1})). Recall
(see, e.g., \cite{daf}) that the system (\ref{pert2}) is called {\it hyperbolic}
if the eigenvalues of the matrix $V^i_j(v)$ are all real and all $n$
eigenvectors are linearly independent. In particular,
{\it strictly hyperbolic} systems are those for which the eigenvalues are all
real and
pairwise distinct.
An
important particular class is the so-called {\it systems of conservation laws}
\beq\label{pert3}
v^i_t + \pal_x \phi^i(v)=0, \quad i=1, \dots, n
\eeq
where the dependent variables are chosen to be densities of conserved
quantities;
the functions $\phi^i(v)$ are the corresponding densities of fluxes (see, e.g.,
\cite{daf} regarding the physical applications of such systems).
The
relationships between solutions of the perturbed and unperturbed systems have
been extensively studied for the case of {\it dissipative} perturbations of
spatially
one-dimensional systems of conservation laws (see, e.g., \cite{bressan} and
the references therein).
Our strategic goal is the study of {\it Hamiltonian} perturbations
of hyperbolic PDEs. Although many concrete examples of such perturbations have
been studied (see, e.g., \cite{whitham, gp, ffm, dm, ll, dn89, kri}), the general concepts and results are still missing.

Let us first explain how to recognize Hamiltonian systems
among all systems of conservation laws. Recall \cite{daf} that the system of
conservation laws (\ref{pert3}) is {\it symmetrizable} in the sense of
Friedrichs and Lax, Godunov, if there exists a
constant symmetric positive definite matrix $\eta=(\eta_{ij})$ such that the matrix
$$
\eta_{is}\frac{\pal \phi^s}{\pal v^j}
$$
is symmetric,
\beq\label{sym}
\eta_{is}\frac{\pal \phi^s}{\pal v^j} = \eta_{js}\frac{\pal \phi^s}{\pal v^i}.
\eeq
In case the symmetry (\ref{sym}) holds true but the symmetric matrix $\eta$ is only 
nondegenerate but
not necessarily positive definite one obtains
{\it weakly symmetrizable} systems of conservation laws.

\begin{lem} The system of conservation laws (\ref{pert3}) is Hamiltonian if it
is weakly symmetrizable.
\end{lem}

\pf Choosing the Poisson brackets in the constant form
$$
\{ v^i(x), v^j(y) \} =\eta^{ij}\delta'(x-y), \quad \left( \eta^{ij}\right) =
\left( \eta_{ij}\right)^{-1}
$$
and a local Hamiltonian
$$
H=\int h(v(x))\, dx
$$
one obtains the Hamiltonian system in the form (\ref{pert3}) with
$$
\phi^i(v) =\eta^{is}\frac{\pal h(v)}{\pal v^s}.
$$
Both sides of (\ref{sym}) coincide then with the Hessian of the Hamiltonian
density $h(v)$. The Lemma is proved. \epf

Recall that weakly symmetrizable systems of conservation laws
enjoy the following important property: they possess two additional conservation
laws
\beq\label{mom}
\pal_t p(v) +\pal_x q(v) =0, \quad p=\frac12 \eta_{ij} v^i v^j, \quad q= v^i
{\pal h\over \pal v^i} - h(v)
\eeq
\beq\label{ene} \pal_t h(v) +\pal_x f(v) =0, \quad f(v) =\frac12 \eta^{ij} 
{\pal h\over \pal v^i}{\pal h\over \pal v^j}
\eeq
where $h(v)$ is the Hamiltonian density in the formulae above. For symmetrizable
systems the function $p(v)$ is nonnegative.

\medskip

The class of Hamiltonian perturbations to be investigated will be written in the form
\beq\label{pert4}
w^i_t +\left\{ w^i(x), H\right\} =w^i_t+ V^i_j(w) w^j_x + \sum_{k\geq 1} \ve^k U^i_k(w; w_x, \dots, w^{(k+1)})=0, \quad i=1, \dots, n
\eeq
where $\ve$ is the small parameter,
$
U^i_k(w; w_x, \dots, w^{(k+1)})
$
are graded homogeneous polynomials\footnote{A different class of perturbations,
for the particular case of the KdV equation, was considered by Y.Kodama
\cite{kod}. In his theory the terms of the perturbative expansion are
polynomials also in $w$ (here $n=1$). The degree on the algebra of differential
polynomials is defined by $\deg u^{(m)}=m+2$, $m\geq 0$. Also some nonlocal terms
appear in the Kodama's perturbation theory. Further developments of this method
can be found in \cite{km}.} in the jet variables $w_x=(w^1_x, \dots,
w^n_x)$, $w_{xx}=(w^1_{xx}, \dots, w^n_{xx})$, \dots, $w^{(k+1)} =\left(w^{1,k+1},
\dots, w^{n,k+1}\right)$
\beq\label{grade}
\deg w^{i,m}=m,\quad i=1,\dots,n,\ m > 0.
\eeq
They arise, e.g., in the study of solutions slowly varying in the space-time
directions \cite{whitham}.  The Hamiltonian are local functionals
\eqa\label{bih2}
&&
H=\int  \left[ h^{[0]}(w) + \ve\, h^{[1]}(w; w_x) + \dots \right] \, dx,
\\
&&
\deg h^{[k]} (w; w_x, \dots, w^{(k)}) =k.
\nn
\eeqa
The Poisson brackets are assumed to be {\it local} in every order in $\ve$, i.e. they are represented as follows
\beq
\{w^i(x),w^j(y)\}=\sum_{m\ge 0}\sum_{l=0}^{m+1} \ve^m
A^{ij}_{m,l}(w;w_x,\dots,w^{(m+1-l)}) \delta^{(l)}(x-y)
 \label{eq-20}
\eeq
with polynomial in the derivatives coefficients
\beq\label{dega0}
\deg A^{ij}_{m,l}(w;w_x,\dots,w^{(m+1-l)})= m-l+1
\eeq
We also assume that the coefficients of
these differential polynomials are ${\mathcal C}^\infty$ smooth functions
on a $n$-dimensional ball $w\in B\subset {\mathbb R}^n$.
It is understood that the antisymmetry and the Jacobi identity for (\ref{eq-20})
hold true as identities for formal power series in $\ve$.
It can be
readily seen that, for an arbitrary local Hamiltonian of the form (\ref{bih2})
the evolutionary systems (\ref{pert4}) has the needed form.

The leading term
\beq\label{pb0}
\left\{ w^i(x), w^j(y)\right\}^{[0]}:= A^{ij}_{0,1}(w(x))\delta'(x-y) +A^{ij}_{0,0}(w(x); w_x(x))\delta(x-y)
\eeq
is itself a  Poisson bracket (the so-called {\it Poisson bracket of hydrodynamic type}, see \cite{dn83}).
We will always assume that
\beq\label{det0}
\det A^{ij}_{0,1}(w)\neq 0
\eeq
for all $w\in B\subset {\mathbb R}^n$. Redenote the coefficients of $\left\{ w^i(x), w^j(y)\right\}^{[0]}$ as follows
\beq\label{metric}
g^{ij}(w):=  A^{ij}_{0,1}(w), \quad Q^{ij}_k(w) w^k_x := A^{ij}_{0,0}(w; w_x)
\eeq
(see (\ref{dega0}).
The coefficient $(g^{ij}(w))$ can be considered as a symmetric nondegenerate bilinear form on the cotangent spaces. The inverse matrix defines a metric
 \beq\label{metric1}
 ds^2 = g_{ij}(w) dw^i dw^j, \quad \left(g_{ij}(w)\right) := \left( g^{ij}(w)\right)^{-1}
 \eeq
 (not necessarily
positive definite).  Recall \cite{dn83} that (\ref{pb0}) - (\ref{metric}) defines a Poisson structure if and only if the metric
is flat and
$Q^{ij}_{k}(w)=-\sum_{l=1}^n g^{il}(w) \Gamma^j_{kl}(w)$ where $\Gamma^j_{kl}$ are the Christoffel symbols of the
Levi-Civita connection of the metric (\ref{metric1}).

The class of the Hamiltonians (\ref{bih2}), Poisson brackets (\ref{eq-20}) and the evolutionary systems (\ref{pert4}) is invariant with respect to {\it Miura-type transformations}
\eqa\label{miura}
&&
w^i\mapsto \tilde w^i = \Phi^i_0(w) + \sum_{k\geq 1} \ve^k \Phi^i_k(w; w_x, \dots, w^{(k)}), \quad i=1, \dots, n
\\
&& \deg \Phi^i_k(w; w_x, \dots, w^{(k)})=k, \quad\det \left( \frac{\pal
\Phi^i_0}{\pal w^j}\right)\neq 0. \nn \eeqa As usual, the
coefficients $ \Phi^i_k(w; w_x, \dots, w^{(k)})$ are assumed to
depend polynomially on the derivatives. Two Poisson brackets of
the form (\ref{eq-20}) are called {\it equivalent} if they are
related by a Miura-type transformation. 

The group of Miura-type transformations is a natural extension of the group of
local diffeomorphisms that plays an important role in the geometrical study of
hyperbolic systems (see, e.g.,\cite{rj}).

An important result of \cite{get} (see also \cite{magri2, DZ1}) says that any Poisson bracket of the form (\ref{eq-20}) can be locally reduced by Miura-type transformations to the constant form
\beq\label{cons}
\left\{ \tilde w^i(x), \tilde w^j(y)\right\} = \eta^{ij}\delta'(x-y), \quad \eta^{ij}=\mbox{const}.
\eeq
We will denote the inverse matrix by the same symbol with lower indices
\beq\label{eta-inv}
\left(\eta_{ij}\right) := \left(\eta^{ij}\right)^{-1}.
\eeq

Connections of the theory of Hamiltonian systems (\ref{pert4}) to the theory
of systems of conservation laws is clear from the following statement.

\begin{lem} \label{lem-scl} By a change of dependent variables of the form (\ref{miura})
the Hamiltonian system (\ref{pert4}) can be recast into the form of a system of conservation laws
\eqa\label{scl}
&&
\tilde w^i_t +\pal_x \psi^i(\tilde w;\tilde w_x, \dots; \ve)=0, \quad i=1, \dots, n
\\
&&
\psi^i(\tilde w; \tilde w_x, \dots; \ve)=\sum_{k\geq 0} \ve^k \psi^i_k(\tilde w; \tilde w_x, \dots, \tilde w^{(k)}),
\nn\\
&&
\deg \psi^i_k(\tilde w; \tilde w_x, \dots, \tilde w^{(k)})=k.
\nn
\eeqa
The system of conservation laws (\ref{scl}) is Hamiltonian with respect to the Poisson bracket
(\ref{cons}) {\em iff} the right hand sides
\beq\label{psi}
\psi_i:= \eta_{ij} \psi^j(\tilde w; \tilde w_x, \dots; \ve)
\eeq
 satisfy
\beq\label{one-main}
\frac{\pal \psi_i} {\pal \tilde w^{j,s}} =\sum_{t\geq s} (-1)^t
\left(\begin{array}{c}  t \\ s\end{array}\right)
\pal_x^{t-s} \frac{\pal \psi_j}{\pal \tilde w^{i,t}}
\eeq
for any $i,\, j=1, \dots, n$, $s=0,\, 1, \, \dots$.
\end{lem}

In this paper we will investigate the structure of
{\it formal} perturbative expansions of the solutions to (\ref{pert4})
\beq\label{pert5}
w^i(x,t;\ve) =v^i(x,t) +\ve\, \delta_1 v^i(x,t) +\ve^2 \delta_2 v^i(x,t) +\dots
\eeq
The leading term solves (\ref{pert2}); the coefficients of the expansion
$\delta_k v^i(x,t)$ are to be determined from linear PDEs with coefficients
depending on $v^i$, $\delta_1 v^i$, \dots, $\delta_{k-1} v^i$ and their
derivatives. Instead of developing this classical technique we propose a
different approach that conceptually goes back to the Poincar\'e's treatment of
perturbative expansions in the celestial mechanics. We will look for a
transformation of the form
\beq\label{pert6}
w^i = v^i +\sum_{k\geq 1} \ve^k \Phi_k^i(v; v_x, \dots, v^{(m_k)}), \quad i=1,
\dots, n
\eeq
that maps {\it any} generic\footnote{We will later be more specific in describing the range of applicability of the transformations (\ref{pert6}).} solution $v^i(x,t)$ of the unperturbed system (\ref{pert2})
to a solution $w^i(x,t;\ve)$ of the perturbed system. An important feature of
such an approach to the perturbation theory is {\it locality}: changing the
functions $v(x,t)$ for a given $t$ only within a small neighborhood of the given
point $x=x_0$ will keep unchanged the values of $w(x,t;\ve)$ outside this
neighborhood. We call (\ref{pert6}) {\it the reducing transformation} for the
perturbed system (\ref{pert4}).

Clearly, applying to (\ref{pert2}) any transformation (\ref{pert6}) {\it polynomial} in the
derivatives (in every order in $\ve$; in that case $m_k=k$) one obtains a
perturbed system of the form (\ref{pert4}). This is the case of {\it trivial}
perturbations.

It is clear that solutions of trivial Hamiltonian perturbations share many properties of the solutions to the unperturbed hyperbolic PDEs (\ref{pert2}). In particular, the trivial perturbation cannot balance the nonlinear effects in the hyperbolic system that typically cause gradient catastrophe of the solution.

\begin{dfn} The system of PDEs (\ref{pert4}) is called {\em quasitrivial} if it is
not trivial but there exists a reducing transformation (\ref{pert6}) with
functions $\Phi^i_k(v; v_x, \dots, v^{(m_k)})$ depending {\em rationally} on the
jet coordinates
\eqa\label{degr}
&&
w^{i,l}, \quad i=1, \dots, n, \quad 1\leq l \leq m_k
\nn\\
&&
\deg \Phi^i_k (v; v_x, \dots, v^{(m_k)}) = k, \quad k\geq 1.
\eeqa
\end{dfn}

The first example of such a reducing transformation was found in \cite{bgi} (see also
\cite{DZ1})
for the KdV equation
\beq\label{kdv}
w_t + w\, w_x + \frac{\ve^2}{12}\, w_{xxx}=0
\eeq
(here $n=1$):
\begin{equation}
w=v +\frac{\epsilon^2}{24} \p_x^2 \, \log v_x
+\epsilon^4 \p_x^2 \left( \frac{v^{(4)}}{1152\, v_x^2}
- \frac{7\, v_{xx} v_{xxx}}{1920\, v_x^3}
+\frac{v_{xx}^3}{ 360\, v_x^4}\right) + O(\epsilon^6).
\end{equation}
It is not an easy
task to check cancellation of all the denominators even in this example!
Because of the denominators the reducing deformation is defined only on the monotone
solutions.

One of the main
outputs of our paper is in proving quasitriviality of  a large class of Hamiltonian perturbations of hyperbolic systems of conservation laws.

The systems in question are {\it bihamiltonian systems of PDEs}. That means
that they can be represented in the Hamiltonian form in two different ways
\beq\label{bih1}
w^i_t = \left\{ H_1,w^i (x)\right\}_1 =\left\{ H_2, w^i(x)\right\}_2, \quad i=1,
\dots, n
\eeq
with two local Hamiltonians (see (\ref{bih2}) above)
and local {\it compatible} Poisson brackets $\{ ~, ~\}_1$, $\{~, ~\}_2$ of the form (\ref{eq-20}). Compatibility means that any linear combination
$$
a_1  \{ ~, ~\}_1+ a_2 \{ ~, ~\}_2
$$
with arbitrary constant coefficients $a_1$, $a_2$ must be again a Poisson
bracket.

The study of bihamiltonian structures was initiated by F. Magri \cite{magri}
in his analysis of the so-called Lenard scheme of constructing the KdV integrals.
I. Dorfman and I. Gelfand \cite{gd} and also A. Fokas and B. Fuchssteiner \cite{FF}
discovered the connections between the bihamiltonian scheme and the theory of
hereditary symmetries of integrable equations. However, it is not easy to apply
these beautiful and
simple ideas to the study of general bihamiltonian PDEs (see the discussion of
the problems occured in \cite{DZ1}).

In this paper we will use a different approach to the study of bihamiltonian
PDEs proposed in \cite{DZ1}. It is based on the careful study of the
transformation properties of the bihamiltonian structures under the
transformations of the form (\ref{pert6}). Let us now proceed to the precise
definitions and
formulations of the results.

We will study bihamiltonian structures defined by compatible pairs of {\it
local} Poisson brackets written in the form of $\ve$-expansion
\begin{eqnarray}
&&\{w^i(x),w^j(y)\}_{a}=\{w^i(x),w^j(y)\}^{[0]}_a
\nn\\&&+\sum_{m\ge 1}\sum_{l=0}^{m+1} \ve^m
A^{ij}_{m,l;a}(w;w_x,\dots,w^{(m+1-l)}) \delta^{(l)}(x-y),
\  a=1,2,\label{eq-2}
\end{eqnarray}
with polynomial in the derivatives coefficients
\beq\label{dega}
\deg A^{ij}_{m,l;a}(w;w_x,\dots,w^{(m+1-l)})= m-l+1.
\eeq
The coefficients of
these differential polynomials are ${\mathcal C}^\infty$ smooth functions
on a $n$-dimensional ball $w\in B\subset {\mathbb R}^n$. Equivalence of bihamiltonian structures is defined with respect to Miura-type transformations.

The leading terms of the bihamiltonian structure is itself a bihamiltonian
structure of the form
\eqa
&&\{w^i(x),w^j(y)\}^{[0]}_a=g^{ij}_a(w(x))\delta'(x-y)+\sum_{k=1}^n Q^{ij}_{a;k}(w(x)) w^k_x
\delta(x-y),\label{eq-1}\\
&&\quad \det(g^{kl}_a)\ne 0 \quad {\mbox{for generic points $w\in M$}},\quad i,j=1,\dots, n,\ a=1,2\nn
\eeqa
(the bihamiltonian structure of the hydrodynamic type).
We additionally assume that
$\det(a_1 g^{kl}_1(w)+a_2 g^{kl}_2(w))$
does not vanish identically for $w\in M$ unless $a_1=a_2=0$.\begin{dfn}
The bihamiltonian structure (\ref{eq-2}) is called {\em semisimple} if the characteristic
polynomial $\det(g^{ij}_2(w)-\lm g^{ij}_1(w))$
in $\lm$ has $n$ pairwise distinct real\footnote{One can relax the requirement
of reality of the roots working with complex manifolds. In that case the
coefficients must be analytic in $w$.}
roots $\lm_1(w),\dots,\lm_n(w)$
for any $w\in B$.
\end{dfn}

The role of semisimplicity assumption can be illustrated by the following

\begin{lem}\label{lem-fs} Given a semisimple bihamiltonian structure
$\{ ~, ~\}^{[0]}_{1,2}$ satisfying the above conditions, denote
$$
\lambda =u^1(w), \dots, \lambda=u^n(w)
$$
the roots of the characteristic equation
\beq\label{roots}
\det\left( g^{ij}_2 (w)-\lambda\, g^{ij}_1(w)\right)=0.
\eeq
The functions $u^1(w)$, \dots, $u^n(w)$ satisfy
$$
\det\left(\frac{\pal u^i(w)}{\pal w^j}\right) \neq 0.
$$
Using these functions as new local coordinates
\beq\label{redc}
w^i=q^i(u^1,\dots, u^n), \quad i=1, \dots, n, \quad \det\left(\frac{\pal q^i(u)}{
\pal u^j}\right)\neq 0
\eeq
reduces both the two
flat metrics to the diagonal form
\beq\label{eq-8}
\frac{\pal u^i}{\pal w^k} \frac{\pal u^j}{\pal w^l}g^{kl}_1(w)=\delta_{ij} f^i(u),\quad
\frac{\pal u^i}{\pal w^k} \frac{\pal u^j}{\pal w^l}g^{kl}_2(w)=\delta_{ij} g^i(u)=\delta_{ij} u^i f^i(u).
\eeq
The coefficients $Q^{ij}_{a;k}$ in the coordinates $(u^1, \dots, u^n)$ read
\eqa
&&\sum_{k=1}^n Q^{ij}_{1;k} u^k_x=\frac12\,\delta_{ij} \p_x f^i+A^{ij},\quad
\sum_{k=1}^n Q^{ij}_{2;k} u^k_x=\frac12\,\delta_{ij} \p_x g^i+B^{ij},\label{eq-9a}\\
&&A^{ij}=\frac12\left(\frac{f^i}{f^j}f^j_iu^{j}_{x}-
\frac{f^j}{f^i}f^i_ju^{i}_{x}\right),\quad
B^{ij}=\frac12\left(\frac{u^i f^i}{f^j}f^j_iu^{j}_{x}-
\frac{u^j f^j}{f^i}f^i_ju^{i}_{x}\right)\label{eq-9}
\eeqa
where $f^i_k=\frac{\p f^i}{\p u^k}$. The leading term of any bihamiltonian system
becomes diagonal in the coordinates $u^1$, \dots, $u^n$,
\beq\label{eq-diag}
u^i_t + V^i(u) u^i_x +{\mathcal O}(\ve)=0, \quad i=1, \dots, n.
\eeq
\end{lem}

Such coordinates are called {\it the canonical coordinates}
of the semisimple bihamiltonian structure.
They are defined up to a permutation. The functions $f_1(u)$, \dots, $f_n(u)$
satisfy a complicated system of nonlinear differential equations. The general
solution to this system depends on $n^2$ arbitrary functions of one variable.
Integrability of this system has recently been proved in \cite{FEV},
\cite{mokhov}. For convenience of the reader we give a brief account of these
results, following \cite{FEV}, in the Appendix below.

The following statement gives a
simple criterion of a Hamiltonian system of conservation laws to be
bihamiltonian.

\begin{lem} \label{lem-fs1} Let us consider a strictly hyperbolic system (\ref{pert2})
Hamiltonian with respect to the Poisson bracket $\{~,~\}_1^{[0]} $. This system is
bihamiltonian with respect to the semisimple Poisson pencil $\{~,~\}^{[0]}_{1,2}$
{\em iff} it becomes diagonal in the canonical coordinates for the
Poisson pencil:
\beq\label{dia}
\frac{\pal u^i}{\pal v^k} \frac{\pal v^l}{\pal u^j} V^k_l(v) = V^i(u) \delta^i_j.
\eeq
\end{lem}

Observe that the canonical coordinates are {\it Riemann invariants} (see, e.g.,
\cite{whitham}) for the leading term of the system of PDEs (\ref{eq-diag}).
The coefficients $V^i(u)$ in the gas dynamics are called {\it characteristic
velocities} \cite{whitham}. In particular the semisimplicity assumption implies
hyperbolicity of the leading term of the Hamiltonian systems.

\begin{dfn}[\cite{DZ1}]
The bihamiltonian structure (\ref{eq-2}) is said to be {\em trivial}
if it can be obtained from the leading term
\beq\label{lead}
\{v^i(x), v^j(y)\}_a^{[0]} = g^{ij}_a(v(x)) \delta'(x-y)
+Q^{ij}_{a;k}(v(x))\delta(x-y), \quad a=1,\, 2
\eeq
by a Miura-type transformation
\eqa\label{triv}
&&
w^i=v^i+\sum_{k\geq 1} \ve^k F^i_k (v; v_x, \dots, v^{(k)}),
\\
&&
\deg F^i_k (v; v_x, \dots, v^{(k)}) = k, \quad i=1, \dots, n
\nn
\eeqa
where the coefficients $F^i_k(v; v_x, \dots, v^{(k)})$ are graded homogeneous
polynomials in the derivatives. It is called
{\em quasitrivial}
if it is not trivial and there exists a transformation
\beq\label{eq-3}
w^i=v^i+\sum_{k\ge 1} \ve^k F^i_k(v;v_x,\dots,v^{(m_k)})
\eeq
reducing (\ref{eq-2}) to (\ref{lead}) but
the functions
$F_k$ depend {\em rationally} on the jet coordinates $v^{i,m},\ m\ge 1$ with
\beq
\deg F_k=k,\quad k\ge 1
\eeq
and  $m_k$ are some positive integers. If such a transformation (\ref{triv}) or (\ref{eq-3})  exists,
it is called a {\em reducing transformation}
of the bihamiltonian structure (\ref{eq-2}).
\end{dfn}

A transformation of the form (\ref{eq-3}) is called a {\em quasi-Miura transformation}.

We are now in a position to formulate the main result of the present paper.

\noindent{\bf Quasitriviality Theorem}
{\em For any semisimple bihamiltonian structure (\ref{eq-2}) there exists a reducing transformation of the form (\ref{eq-3}). The coefficients $F^i_k$ have the form
\eqa\label{3g-2}
&&
 F^i_k(v;v_x,\dots,v^{(m_k)})\in C^\infty(B)\left[v_x, \dots, v^{(m_k)}\right]\,
 \left[ \left( u^1_x u^2_x \dots u^n_x\right)^{-1}\right]
 \\
 &&
 m_k\leq \left[ \frac{3\, k}{2}\right].
\nn
\eeqa
Here $u^i=u^i(v)$ are the canonical coordinates (\ref{eq-8}).
}
\vskip 0.2cm
Using this theorem we achieve the goal of constructing the reducing
transformation for a bihamiltonian system (\ref{bih1}), (\ref{bih2}):

\begin{cor}\label{thm2}
The reducing transformation for the bihamiltonian structure (\ref{eq-2})
is also a reducing transformation for any
bihamiltonian system (\ref{bih1}).
\end{cor}

Another corollary says that the solution of any system of bihamiltonian PDEs
of the above form can be reduced to solving {\em linear PDEs}. Let us first
rewrite the reducing transformation in the canonical coordinates
\beq\label{red-can}
\tilde u^i= u^i +\sum_{k\geq 1} \ve^k G^i_k(u; u_x, \dots, u^{(m_k)}), \quad i=1,
\dots, n.
\eeq

Let $W^i(u; \ve)$, $i=1, \dots, n$ be an arbitrary solution to the linear
system
\beq\label{tsar}
\frac{\pal W^i}{\pal u^j} =\frac{\pal V^i /\pal u^j}{ V^i - V^j}
\left( W^i-W^j\right),
\quad i\neq j
\eeq
in the class of formal power series in $\ve$. In this system the functions
$V^i(u)$ are eigenvalues of the matrix $V^i_j(w)$, cf. (\ref{eq-diag}).
Let us assume that the system
of equations
\beq\label{hod}
x= V^i(u)t +W^i(u;\ve=0), \quad i=1, \dots, n
\eeq
has a solution $(x,t,u)=(x_0, t_0, u_0)$  such that
\beq\label{jac}
\det\left( t\,\pal V^i(u) /\pal u^j + \pal W^i(u;\ve)/\pal u^j\right)_{t=t_0, \,
u=u_0,\, \ve=0} \neq 0.
\eeq
For $(x,t)$ sufficiently closed to $(x_0, t_0)$ denote
$u(x,t)=(u^1(x,t), \dots, u^n(x,t))$ the unique solution to the equations
(\ref{hod}) such that
$$
u(x_0, t_0)=u_0.
$$
Applying the transformation (\ref{red-can}) to the vector-function $u(x,t)$
we obtain a vector-function $\tilde u(x,t;\ve)$. Finally the substitution
\beq\label{resh}
q^i(\tilde u(x,t;\ve))=: w^i(x,t;\ve), \quad i=1, \dots
\eeq
yields $n$ functions $w^1(x,t;\ve)$, \dots, $w^n(x,t;\ve)$. Here the functions
$q^i(u)$ are defined as in (\ref{redc}).

\begin{cor} \label{cor-tsar} The functions (\ref{resh}) satisfies (\ref{bih1}). Conversely, any {\em monotone} at $x=x_0, ~t=t_0$ solution to
(\ref{bih1}) can be obtained by this procedure.
\end{cor}

By definition, the solution $w(x,t;\ve)$ is called monotone if all the $x$-derivatives
$$
\pal_xu^1(w(x,t;\ve)), \dots , \pal_xu^n(w(x,t;\ve))
$$
do not vanish for $x=x_0$, $t=t_0$, $\ve=0$.

Finally, we can combine the Quasitriviality Theorem  with the main result of the
recent paper \cite{LZ1} in order to describe the complete set of invariants of
bihamiltonian structures of the above form with the {\it given} leading term
$\{w^i(x), w^j(y)\}^{[0]}_a$.

Introduce the following combinations of the coefficients of $\ve
\,\delta''(x-y)$ and $\ve^2\delta'''(x-y)$ of the bihamiltonian structure
\beq\label{ss-31}
P^{ij}_a(u)=\frac{\pal u^i}{\pal w^k} \frac{\pal u^j}{\pal w^l}A^{kl}_{1,2;a}(w),
\ Q^{ij}_a(u)=\frac{\pal u^i}{\pal w^k} \frac{\pal u^j}{\pal w^l} A^{kl}_{2,3;a}(w),\ i,j=1,\dots,n,
a=1,2.
\eeq
Define the functions
\beq\label{fc}
c_i(u)=\frac1{3 (f^i(u))^2} \left(Q^{ii}_2- u^i Q^{ii}_1+\sum_{k\ne i}\frac{(P^{ki}_2-u^i P^{ki}_1)^2}{f^k(u)
(u^k-u^i)}\right),\quad i=1,\dots,n.
\eeq
The functions $c^i(u)$ are called {\it central invariants} of the bihamiltonian
structure (\ref{eq-2}).
The main result of \cite{LZ1} on the classification of infinitesimal deformations of bihamiltonian
structures of hydrodynamic type can be reformulated as follows
\begin{cor}\label{thm51} i) Each function $c_i(u)$ defined in (\ref{fc}) depends only on $u^i$.\quad
ii) Two semisimple bihamiltonian structures (\ref{eq-2}) with the same leading terms $\{\ ,\ \}_a^{[0]},
\ a=1,2$ are equivalent {\em iff} they have the same set of central invariants $c_i(u^i),  i=1,\dots,n$. In particular, a polynomial in the derivatives reducing transformation exists
{\em iff} all the central invariants vanish:
$$
c_1=c_2=\dots=c_n=0.
$$
\end{cor}

The papers is organized as follows.
In Section \ref{sec-2-a} we recall some basic notions of the theory of 
Poisson structures for PDEs and prove Lemma \ref{lem-scl}, \ref{lem-fs}
and Lemma \ref{lem-fs1}.
In Section \ref{sec-2} and Section \ref{sec-3} we give the proofs of the
Quasitriviality
Theorem and the Corollaries \ref{thm2}, \ref{cor-tsar}. In Section \ref{sec-4} 
we 
reformulate the main result of \cite{LZ1}
on the classification of infinitesimal deformations of a semisimple bihamiltonian structure of
hydrodynamic type and prove Corollary \ref{thm51}.  In the final section we give some examples of bihamiltonian structures of the
class studied in this paper and formulate some open problems. In the Appendix we
briefly present, following \cite{FEV}, the theory of semisimple bihamiltonian
structures of hydrodynamic type.
\vskip 0.5truecm
\noindent{\bf Acknowledgments.} The researches of B.D. were
partially supported by European Science Foundation Programme ``Methods of
Integrable Systems, Geometry, Applied Mathematics" (MISGAM). The researches of Y.Z. were partially supported
by the Chinese National Science Fund for Distinguished Young Scholars grant
No.10025101 and the Special Funds of Chinese Major Basic Research Project
``Nonlinear Sciences''. Y.Z. and S.L. thank
Abdus Salam International Centre for Theoretical Physics and SISSA
where part of their work was done for the hospitality.

\section{Some basic notions about Poisson structures for PDEs}\label{sec-2-a}
Like in finite dimensional Poisson geometry, an infinite dimensional Poisson structure of the form
(\ref{eq-1}) or (\ref{eq-2}) can be represented by a local bivector on the formal loop space
of the manifold $M$. Recall that in our considerations the manifold $M$ will always be a $n$-dimensional ball.
In general, let $w^1,\dots,w^n$ be a local coordinate system of a chart of the
manifold $M$. A {\it local translation invariant} $k$-vector \cite{DZ1} is a formal infinite sum of the form
\begin{equation}
\al=\sum \frac1{k!}\pal_{x_1}^{s_1}\dots\pal_{x_k}^{s_k} A^{i_1\dots i_k}
\frac{\pal}{\pal w^{i_1,s_1}(x_1)}\wedge\dots \wedge\frac{\pal}{\pal w^{i_k,s_k}(x_k)}.
\end{equation}
Here the coefficients $A$'s have the expressions
\begin{equation}
A^{i_1\dots i_k}=\sum_{p_2,\dots,p_{k}\ge 0} B^{i_1\dots i_k}_{p_2\dots p_k}(w(x_1);
w_x(x_1),\dots) \delta^{(p_2)}(x_1-x_2)\dots \delta^{(p_k)}(x_1-x_k).
\end{equation}
with only a finite number of nonzero terms in the summation,
and $B^{i_1\dots i_k}_{p_2\dots p_k}(w; w_x, \dots)$
is a smooth function on a domain in the jet space
$J^N(M)$ for certain integer $N$ that depends on the indices $i_1, \dots, i_k$ and $p_2, \dots, p_k$.
The delta-function and its derivatives are defined formally by
\beq\label{del1}
\int f(w(y); w_y(y), w_{yy}(y), \dots) \, \delta^{(k)}(x-y)\, dy =
\pal_x^k f(w(x); w_x(x), w_{xx}(x), \dots).
\eeq
In this formula the operator of total derivative $\pal_x$ is defined by
\beq\label{palx}
\pal_x f(w; w_x, w_{xx}, \dots) =\sum w^{i,s+1} \frac{\pal f}{\pal w^{i,s}}.
\eeq
Note the useful identity
\eqa\label{del2}
&&
f(w(y); w_y(y), w_{yy}(y), \dots) \, \delta^{(k)}(x-y)
\nn\\
&&
=\sum_{m=0}^k \left(\begin{array}{c}k \\ m\end{array}\right)\, \pal_x^{m}
f(w(x); w_x(x), w_{xx}(x), \dots)\, \delta^{(k-m)}(x-y).
\eeqa

The distributions
\begin{equation}
A^{i_1\dots i_k}=A^{i_1\dots i_k}(x_1,\dots,x_k;w(x_1),\dots,w(x_k),\dots)
\end{equation}
are antisymmetric with respect to the simultaneous permutations
$
i_p,x_p\leftrightarrow i_q, x_q.
$
They are called {\it the components} of the local $k$-vector $\al$.
Note that in the definition of local
$k$-vectors given in \cite{DZ1} it is required that the functions $B^{i_1\dots i_k}_{p_2\dots p_k}$
are differential polynomials. Here we drop this requirement for the convenience of our use of these notations
during our proof of the theorem.
The space of all such local $k$-vectors is still denoted by $\Lambda_{loc}^k$ as it is done in \cite{DZ1}.
For $k=0$ by definition $\Lambda^0_{loc}$ is the space of local functionals of the form
\beq
I=\int f(w,w_x,\dots,w^{(m)})\, dx
\eeq
The Schouten-Nijenhuis bracket is defined on the space of local multi-vectors
\beq
[\,,\,]:\ \Lambda^k_{loc}\times \Lambda^l_{loc}\to \Lambda^{k+l-1}_{loc},\quad k,l\ge 0.
\eeq
It generalizes the usual commutator of two local vector fields and possesses the following properties
\eqa
&&[\al,\beta]=(-1)^{kl} [\beta,\al],\\
&&(-1)^{km}[[\al,\beta],\gamma]+(-1)^{kl}[[\beta,\gamma],\al]+(-1)^{lm} [[\gamma,\al],\beta]=0
\eeqa
for any $\al\in\Lambda^k_{loc},\, \beta\in\Lambda^l_{loc},\, \gamma\in\Lambda^m_{loc}$.
For the definition of the Schouten-Nijenhuis
bracket see \cite{DZ1} and references therein. Here we write down the formulae,
used below, for the bracket of a local bivector with a local functional
and with a local vector field. Let a local vector field $\xi$
and a local bivector $\varpi$ have the representation
\eqa\label{eq-4}
&&\xi=\sum_{i=1}^n\sum_{s\ge 0} \pal_x^s \xi^{i}(w(x);w_x(x),\dots,w^{(m_i)})\frac{\pal}{\pal w^{i,s}(x)},\\
&&
\varpi=\frac12\sum \pal_{x}^{s} \pal_y^t \varpi^{ij}\frac{\pal}{\pal w^{i,s}(x)}\wedge\frac{\pal}{\pal w^{j,t}(y)}.
\label{eq-5}
\eeqa
Here we assume that
\beq\label{redu}
\varpi^{ij}=\sum_{k\ge 0} A^{ij}_k(w(x),w_x(x),\dots,w^{(m_k)}(x))\, \delta^{(k)}(x-y).
\eeq
Then the components of $[\varpi,I]$ and of $[\varpi, \xi]$ are given respectively by
\eqa
&& [\varpi,I]^i=\sum_{j,k} A^{ij}_k \pal_x^k\frac{\delta I}{\delta w^j(x)},\\
&&[\varpi, \xi]^{ij}=\sum_{k,t}\left(\pal_x^t \xi^k(w(x); \dots)
\frac{\pal \varpi^{ij}}{\pal w^{k,t}(x)}
-\frac{\pal \xi^i(w(x); \dots)}{\pal w^{k,t}(x)} \pal_x^t \varpi^{kj}\right.\nn\\
&&\quad \quad \qquad\left. -\frac{\pal
\xi^j(w(y); \dots)}{\pal w^{k,t}(y)} \pal_y^t \varpi^{ik}\right).
\label{lie-der}
\eeqa
In the last formula it is understood that
$$
\pal_y^t \varpi^{ij}=\sum_{k\ge 0}(-1)^t
A^{ij}_k(w(x),w_x(x),\dots,w^{(m_k)}(x))\, \delta^{(k+t)}(x-y).
$$
and the identity (\ref{del2}) has been used in order to represent the resulting
bivector in the normalized form (\ref{redu}).

Let us denote by $\varpi_1, \varpi_2$ the two bivectors that correspond to the bihamiltonian
structure (\ref{eq-1}), the components $\varpi^{ij}_a, a=1,2$ are given by the right hand side of
(\ref{eq-1}). The bihamiltonian property is equivalent to the following identity
that is valid for an arbitrary parameter $\lm$:
\beq\label{biha}
[\varpi_2-\lm\, \varpi_1,\varpi_2-\lm\, \varpi_1]=0.
\eeq
Denote by $\pal_1, \pal_2$ the differentials associated
with $\varpi_1, \varpi_2$.
By definition
\beq
\pal_a:\ \Lambda^k_{loc}\to \Lambda^{k+1}_{loc},\quad \pal_a \al=[\varpi_a,\al],\quad
\forall \al\in \Lambda^k_{loc},\quad a=1,2.
\eeq
The bihamiltonian property (\ref{biha}) can be recast in the form
\beq\label{biha1}
\pal_1^2 =\pal_2^2=\pal_1\pal_2+\pal_2\pal_1=0.
\eeq
The important fact that we need to use below is vanishing of the
first and second Poisson cohomologies
\beq
H^k({\cal {L}}(M),\varpi_a)=
\left.{\mbox{Ker}}\,\pal_a\right|_{\Lambda^k_{loc}}/\left.{\mbox{Im}}\,\pal_a \right|_{\Lambda^{k-1}_{loc}},
\quad a=1,2,\ k=1,2.
\eeq
This fact is proved in \cite{get, magri2, DZ1}.
It readily implies, along with the results of \cite{dn83}
the reducibility of any Poisson bracket of the form (\ref{eq-20}) - (\ref{det0})
to the constant form (\ref{cons}).

Let us now give

\noindent{\em Proof of  Lemma \ref{lem-scl}}. For the
Poisson bracket written in the form (\ref{cons}) the Hamiltonian system reads
$$
\tilde w^i_t =\{ H, \tilde w^i(x)\} = -\eta^{ij} \pal_x\frac{\delta H}{\delta \tilde
w^j(x)}.
$$
This gives a system of conservation laws form with
$$
\psi^i=\eta^{ij} \frac{\delta H}{\delta \tilde
w^j(x)}.
$$
The equations (\ref{one-main}) is nothing but the spelling of the classical
Volterra criterion \cite{volterra} for functions $\psi_i$ to be representable in the form of
variational derivatives. \epf

We pass now to the theory of canonical coordinates.

\noindent{\em Proof of Lemma \ref{lem-fs}}.
{}From the result of \cite{FEV, mokhov} it follows the existence of a system of
local coordinates $\hat u^1,\dots,\hat u^n$ such that both of the metrics become diagonal and
$
g^{ij}_1=\delta_{ij} h^i(\hat u),\ g^{ij}_2=\delta_{ij}\lm_i(\hat u^i) h^i(\hat u).
$
Since by our assumption
$\det(a_1 g^{kl}_1(w)+a_2 g^{kl}_2(w))$
does not vanish identically for $w\in M$ unless $a_1=a_2=0$, we can choose
$u^1=\lm_1(\hat u^1),\dots,u^n=\lm_n(\hat u^n)$ as a system of local coordinates which are just the canonical
coordinates. \epf

We now proceed to

\noindent{\em Proof of Lemma \ref{lem-fs1}}.
In the canonical coordinates a bihamiltonian system
$$
u^i_t =\{ H_1 , u^i(x)\}_1^{[0]} = \{ H_2 , u^i(x)\}_2^{[0]},
$$
$$
H_a =\int h_a^{[0]}(u)\, dx, \quad a=1, \, 2
$$
has the expression
\beq
u^i_t=-\sum_{j=1}^n V^i_j(u) u^j_x,\quad i=1,\dots,n
\eeq
where
\beq\label{riv}
V^i_j(u)=f^i(u) {\mathcal A}_{ij} h^{[0]}_1(u)=u^i f^i(u) {\mathcal A}_{ij} h^{[0]}_2(u),\quad {\text {for}}\ i\ne j.
\eeq
Here the linear differential operators ${\mathcal A}_{ij}$ are defined by
\beq
{\mathcal A}_{ij}=\frac{\pal^2}{\pal u^i\pal u^j}+\frac12 \frac{\pal\left(\log f^i(u)\right)}
{\pal u^j}\frac{\pal}{\pal u^i}+
\frac12 \frac{\pal\left(\log f^j(u)\right)}{\pal u^i}\frac{\pal}{\pal u^j}.
\eeq
Symmetry with respect to to the indices $i, j$ implies
\beq
(u^i-u^j) {\mathcal A}_{ij} h^{[0]}_2(u)=0,\quad i\ne j.
\eeq
Thus $V^i_j(u)=0$ when $i\ne j$. This proves the first part of the lemma.

To prove the converse statement we use the following result of \cite{tsarev}:
the diagonal system
$$
u^i_t +V^i(u) u^i_x=0, \quad i=1, \dots, n
$$
is Hamiltonian with respect to to the Poisson bracket associated with the diagonal metric
of zero curvature
$$
ds^2 =\sum_{i=1}^n g_{ii}(u) (du^i)^2
$$
{\em iff} the following equations holds true
\beq\label{tsa}
\pal_k V^i(u) = \left(V^k(u) - V^i(u)\right)\, \pal_k \log \sqrt{ g_{ii}(u)},
\quad i\neq k.
\eeq
By assumption these equations hold true for the first metric
$$
g_{ii}(u)= 1/ f_i(u).
$$
For the pencil $\{~,~\}_2^{[0]} -\lambda \, \{~,~\}_1^{[0]}$
one has to replace
$$
g_{ii}(u) = \frac{1}{(u^i-\lambda) f^i(u)}.
$$
Such a replacement does not change the equations (\ref{tsa}). The Lemma is
proved. \epf

\section{Proof of the Quasitriviality Theorem}\label{sec-2}

In this and the forthcoming sections we assume that the bihamiltonian structure $(\varpi_1,\varpi_2)$
defined by (\ref{eq-1})
is semisimple and we
work in canonical coordinates $u^1$,\dots, $u^n$.
We first formulate and prove some intermediate results crucial for the proof
of the Quasitriviality Theorem.
\begin{thm}\label{lem-4}
Assume that a vector fields $X$ has components of the form
\beq
X^i=\sum X^i_j(u)\,u^j_x,\quad i=1,\dots,n
\eeq
and satisfies
\beq\label{d1d2-b}
\pal_1 \pal_2 X=0.
\eeq
Then there exist two local functionals $I, J$ of the form
\beq\label{gtg}
I=\int G(u(x)) dx,\quad J=\int \tilde G(u(x)) dx
\eeq
such that $X$ has the representation $X=\pal_1 I-\pal_2 J$.
\end{thm}

\begin{prf}
In this proof summations over repeated Greek indices are assumed. Let us redenote the components
of the two bivectors that correspond to the bihamiltonian structure (\ref{eq-1})
in the form
\eqa
&&\varpi_1^{ij}=g^{ij}\delta'(x-y)+\Gamma^{ij}_\al u^\al_x \delta(x-y), \nn\\
&&\varpi_2^{ij}=\tg^{ij}\delta'(x-y)+\tGamma^{ij}_\al u^\al_x \delta(x-y),\nn
\eeqa
The Levi-Civita connections of these two metrics $g^{ij}, \tg^{ij}$
are denoted by $\nabla$ and $\tnabla$
respectively.
Denote $\nabla_i,\tnabla_i$ the covariant derivatives
of these two connections along $\frac{\p}{\p u^i}$.
We also introduce the notations
$\nabla^i=g^{i\al}\nabla_\al,\tnabla^i=\tg^{i\al}\tnabla_\al$.

The condition (\ref{d1d2-b}) implies existence of a vector field $Y$ with components of the form
\beq
Y^i=\sum Y^i_j(u) u^j_x,\quad i=1,\dots,n
\eeq
such that $\pal_1 X=\pal_2 Y$.
Denote by
\beq\label{eq-z}
Z^{ij}=\sum_{p\ge 0} Z^{ij}_p(u(x),u_x(x),\dots)\, \delta^{(p)}(x-y):=\left(
\pal_1 X-\pal_2 Y\right)^{ij}
\eeq
the components of the bivector $\pal_1 X-\pal_2  Y $, and by $Z^{ij}_{p,(k,m)}$ the derivatives
$\frac{\pal Z^{ij}_p}{\pal u^{k,m}}$. Then we have
\eqa
Z^{ij}_2&=&(X^{ij}-X^{ji})-(Y^{ij}-Y^{ji})=0, \label{n1cond-1}\\
Z^{ij}_{0,(k,2)}&=&\left(\nabla_k X^{ij}-\nabla^i X^j_k+
\Gamma^j_{k\al}\left(X^{\al i}-X^{i\al}\right)\right)\nn\\
&&\quad-\left(\tnabla_k Y^{ij}-\tnabla^i Y^j_k+
\tGamma^j_{k\al}\left(Y^{\al i}-Y^{i\al}\right)\right)=0 \label{n1cond-2}
\eeqa
where $X^{ij}=g^{i\al}X^j_{\al}, Y^{ij}=\tg^{i\al}Y^j_{\al}$.
{}From the above two equations we obtain
\beq
\tnabla^kY^{ij}-\tnabla^iY^{kj}=\tg^{k\al}\left(\nabla_\al X^{ij}-\nabla^iX^j_\al+
T^j_{\al\bt}(X^{i\bt}-X^{\bt i})\right). \label{n1cond-3}
\eeq
Here the components of the $(1,2)$-tensor $T$ are defined by
$T^j_{\al\bt}=\tGamma^j_{\al\bt}-\Gamma^j_{\al\bt}$.
Since the l.h.s of the above equation is antisymmetric with respect to $k,i$, we have
\eqa
&&\tg^{k\al}\left(\nabla_\al X^{ij}-\nabla^iX^j_\al+T^j_{\al\bt}(X^{i\bt}-X^{\bt i})\right)+\nn\\
&&\quad\tg^{i\al}\left(\nabla_\al X^{kj}-\nabla^kX^j_\al+T^j_{\al\bt}(X^{k\bt}-X^{\bt k})\right)=0.
\label{n1cond-sym}
\eeqa
The following trivial identity
\eqa &&
\left(\tnabla^kY^{ij}-\tnabla^iY^{kj}\right)+
\left(\tnabla^iY^{jk}-\tnabla^jY^{ik}\right)+
\left(\tnabla^jY^{ki}-\tnabla^kY^{ji}\right) \nn\\
&&\quad=
\tnabla^k\left(Y^{ij}-Y^{ji}\right)+
\tnabla^i\left(Y^{jk}-Y^{kj}\right)+
\tnabla^j\left(Y^{ki}-Y^{ik}\right), \nn
\eeqa
implies that
\eqa
&&\tnabla^k\left(X^{ij}-X^{ji}\right)+\tnabla^i\left(X^{jk}-X^{kj}\right)+\tnabla^j\left(X^{ki}-X^{ik}\right) \nn\\
&=&\tg^{k\al}\left(\nabla_\al X^{ij}-\nabla^iX^j_\al+T^j_{\al\bt}(X^{i\bt}-X^{\bt i})\right)+\nn\\
&&\tg^{i\al}\left(\nabla_\al X^{jk}-\nabla^jX^k_\al+T^k_{\al\bt}(X^{j\bt}-X^{\bt j})\right)+\nn\\
&&\tg^{j\al}\left(\nabla_\al X^{ki}-\nabla^kX^i_\al+T^i_{\al\bt}(X^{k\bt}-X^{\bt k})\right).\label{n1cond-4}
\eeqa
By using the formula
\beq
\tnabla_kA^{ij}=\nabla_kA^{ij}+T^i_{k\al}A^{\al j}+T^j_{k\al}A^{i\al} \nn
\eeq
we can simplify the equations (\ref{n1cond-4}) to the form
\eqa
&&\tg^{k\al}\left(\nabla_\al X^{ji}-\nabla^iX^j_\al+T^i_{\al\bt}(X^{j\bt}-X^{\bt j})\right)+\nn\\
&&\tg^{i\al}\left(\nabla_\al X^{kj}-\nabla^jX^k_\al+T^j_{\al\bt}(X^{k\bt}-X^{\bt k})\right)+\nn\\
&&\tg^{j\al}\left(\nabla_\al X^{ik}-\nabla^kX^i_\al+T^k_{\al\bt}(X^{i\bt}-X^{\bt i})\right)=0.\label{n1cond-cyc}
\eeqa

Let us employ the equations (\ref{n1cond-sym}) and (\ref{n1cond-cyc}) to prove the existence two local functionals
$I, J$ of the form (\ref{gtg})
such that $X=\pal_1 I-\pal_2 J$. Equivalently, we need to find functions $G, \tilde G$ that satisfy the conditions
\beq \label{gFgF0}
X^i_j=g_{j\al}\nabla^i\nabla^\al G-\tg_{j\al}{\tilde\nabla}^i{\tilde\nabla}^\al {\tilde G}
\eeq
To this end, we first define two symmetric $(2,0)$-tensors $A, \tilde A$ such that
\beq \label{gFgF}
X^i_j=g_{j\al} A^{\al i}-\tg_{j\al}\tilde{A}^{\al i}
\eeq
In the canonical coordinates the off-diagonal components of $A,\tilde{A}$ are uniquely determined
by the above relations and have the explicit forms
\beq
A^{ij}=\frac{g^iX^j_i-g^jX^i_j}{u^i-u^j},\ \tilde{A}^{ij}=\frac{u^jg^iX^j_i-u^ig^jX^i_j}{u^i-u^j}
\quad \textrm{for}\ i \ne j.\label{zh9}
\eeq
Here we use the fact that in the canonical coordinates the two metrics have components of the form
$g^{ij}=\delta_{ij} f^i,\ {\tilde g}^{ij}=\delta_{ij} g^i=\delta_{ij} u^i f^i$.
For an arbitrary choice of the diagonal components $A^{ii}$ the above relation uniquely determines
the diagonal components ${\tilde A}^{ii}$ by
\beq\label{zh6}
{\tilde A}^{ii}=u^i \left( f^i X^i_i-A^{ii}\right).
\eeq
We will specify the choice of $A^{ii}, i=1,\dots, n$
in a moment. Let us now express the equations (\ref{n1cond-sym}) and (\ref{n1cond-cyc}) in terms of
the components of the tensors $A, \tilde A$. By substituting the the expression (\ref{gFgF}) of $X^i_j$ into
(\ref{n1cond-sym}) and (\ref{n1cond-cyc}) and by using
 the fact that $g_{\al\bt}, \tg_{\al\bt}, T^i_{\al\bt}$ are all diagonal
with respect to $\al,\bt$ in the canonical coordinates, we arrive at
\eqa
&&\left(u^k-u^i\right)\left(\nabla^kA^{ij}-\nabla^iA^{kj}\right)+
\left(\frac1{u^k}-\frac1{u^i}\right)\left(\tnabla^k\tilde{A}^{ij}-\tnabla^i\tilde{A}^{kj}\right)=0,\label{zh2}\\
&&u^k\left(\nabla^kA^{ij}-\nabla^iA^{kj}\right)-
\left(\frac1{u^j}\tnabla^k\tilde{A}^{ij}-\frac1{u^i}\tnabla^i\tilde{A}^{kj}\right)\nn\\
&&\qquad +u^i\left(\nabla^iA^{jk}-\nabla^jA^{ik}\right)-
\left(\frac1{u^k}\tnabla^i\tilde{A}^{jk}-\frac1{u^j}\tnabla^j\tilde{A}^{ik}\right) \nn\\
&&\qquad+ u^j\left(\nabla^jA^{ki}-\nabla^kA^{ji}\right)-
\left(\frac1{u^i}\tnabla^j\tilde{A}^{ki}-\frac1{u^k}\tnabla^k\tilde{A}^{ji}\right)=0.\label{zh1}
\eeqa
Rewrite (\ref{zh1}) into the form
\eqa
&&u^k\left(\nabla^kA^{ij}-\nabla^iA^{kj}\right)+
\frac1{u^k}\left(\tnabla^k\tilde{A}^{ij}-\tnabla^i\tilde{A}^{kj}\right)\nn\\
&&\quad +\textrm{terms obtained by
cyclic permutations of $(i,j,k)$}=0.\label{zh3}
\eeqa
By using (\ref{zh2}), we can replace $u^k$ by $u^i$ in the first two terms of (\ref{zh3}). Then after
the cancellation of some terms we arrive at the simplification of (\ref{zh1})
\beq\label{zh4}
\left(u^i-u^j\right)\left(\nabla^j A^{ki}-\nabla^k A^{ij}\right)+
\left(\frac1{u^i}-\frac1{u^j}\right)\left(\tnabla^j \tilde{A}^{ki}-\tnabla^k\tilde{A}^{ij}\right)=0.
\eeq
Changing the indices $(i,j,k)\mapsto (j,k,i)$ we obtain
\beq
\left(u^j-u^k\right)\left(\nabla^kA^{ij}-\nabla^iA^{kj}\right)+
\left(\frac1{u^j}-\frac1{u^k}\right)\left(\tnabla^k\tilde{A}^{ij}-\tnabla^i\tilde{A}^{kj}\right)=0.\label{zh5}
\eeq
{}From the equations (\ref{zh2}), (\ref{zh5}) it readily follows that
\beq\label{zh8}
\nabla^kA^{ij}=\nabla^iA^{kj},\ \tnabla^k\tilde{A}^{ij}=\tnabla^i\tilde{A}^{kj}\quad {\mbox {for}}\ i\ne j, \
k\ne j.
\eeq
Now let us proceed to choosing the diagonal components $A^{ii}$ in such a way
to ensure that the components of the tensor
$\nabla^kA^{ij}$ are totaslly symmetric in $i$, $j$, $k$. This amounts to require that $A^{ii}$
should satisfy
\beq
\nabla^kA^{ii}=\nabla^i A^{ki}, \ i,k=1,\dots,n,\  k\ne i.\label{zh7}
\eeq
The existence of solutions $A^{ii}$ is guaranteed by the compatibility of the above systems due to the
equalities
\beq
\nabla^j\left(\nabla^i A^{ki}\right)=\nabla^k\left(\nabla^i A^{ji}\right),\quad
{\mbox{for distinct $i,j,k$}}.
\eeq
Fix a solution $A^{ii}, i=1,\dots,n$ of the system (\ref{zh7}). {}From the validity of the equations (\ref{zh2})
and (\ref{zh8})  we
know that the tensor ${\tilde A}$ with components ${\tilde A}^{ij}$ determined by (\ref{zh9}) and (\ref{zh6})
also have the property of symmetry of  $\tnabla^k {\tilde A}^{ij}$ in $i$, $j$,
$k$.
Thus we can find functions $G(u), \tilde G(u)$ such that
$$A^{ij}=\nabla^i\nabla^jG,\ \tilde{A}^{ij}=\tnabla^i\tnabla^j\tilde{G}.$$
The lemma is proved.
\end{prf}

The above theorem implies that the linear in $\ve$ terms of the bihamiltonian structure (\ref{eq-2})
can be eliminated by a Miura-type transformation.

Denote by ${\cal A}$ the space of functions that can be represented as
a finite sum of rational functions of the form
\beq\label{mono}
\frac{P^i_{j_1,\dots,j_m}(u;u_x,\dots)}{u^{j_1}_x\dots u^{j_m}_x}, \quad m\geq 0
\eeq
(no denominator for $m=0$).
Here $P^i_{j_1,\dots,j_m}$ are quasihomogeneous differential polynomials. Define a gradation on the ring ${\cal A}$
by
\beq\label{grade-b}
\deg u^{i,m}=m,\quad i=1,\dots,n,\ m\ge 0.
\eeq
We call elements of ${\cal A}$ {\em almost differential polynomials}. Below we will also
encounter functions that belong to the ring
$$
{\tilde {\cal A}}={\cal A}[\log u^1_x,\dots,\log u^n_x].
$$
It is also a graded
ring with the definition of degree (\ref{grade}) and
\beq\label{grade-log}
\deg( \log u^i_x)=0,\quad i=1,\dots,n.
\eeq

\begin{thm}\label{thm3}
Let $X \in \Lambda^1_{loc}$ be a local vector field with components
\beq
X^i(u,u_x,\dots,u^{(N)}),\ i=1,\dots,n,\ N \ge 1,
\eeq
where $X^i$ are homogeneous almost differential polynomials of degree $d\ge 3$.
If $X$ satisfies the condition
\beq\label{d1d2}
\pal_1 \pal_2 X=0
\eeq
then there exist local functionals
\beq
I=\int f(u,u_x,\dots,u^{([\frac{N}2])}) dx,\quad J=\int g(u,u_x,\dots,u^{([\frac{N}2])}) dx,
\eeq
with densities that are homogeneous almost differential polynomials of degree $d-1$
such that $X=\pal_1 I-\pal_2 J$.
\end{thm}

We will prove the Quasitriviality Theorem
by induction on the highest order of the $x$-derivatives of $u^k$ on which the components $X^i$
of the vector field $X$  depend. The following lemmas spell out some important properties
of the vector field $X$ that are implied by the condition (\ref{d1d2}).

In what follows, for a function $A=A(u,u_x,\dots)$ we will use the subscript $(k,m)$ to indicate
the derivative of $A$ with respect to $u^{k,m}$, i.e., $A_{(k,m)}=\frac{\pal A}{\pal u^{k,m}}$.
\begin{lem} \label{lem-1}
For any two vector fields $X,  Y $ with components of the form
\beq\label{eq-11}
X^i=X^i(u,u_x,\dots,u^{(N)}),\quad  Y ^i=Y^i(u,u_x,\dots,u^{(N)}),\quad N\ge 1,
\eeq
vanishing of $Z^{ij}_{0,(k,2N+1)}$ and $Z^{ij}_{0,(k,2N)}$, where $Z^{ij}$ are
defined as in (\ref{eq-z}),
implies that the components of the vector fields $X, Y$ must take the form
\eqa
&& X^i=\sum_{j=1}^n \left(u^j G^i_j(u,\dots,u^{(N-1)})+F^i_j(u,\dots,u^{(N-1)})\right)u^{j,N}
+Q^i(u,\dots,u^{(N-1)}),\nn\\
&& Y^i=\sum_{j=1}^n G^i_j(u,\dots,u^{(N-1)})\,u^{j,N}
+R^i(u,\dots,u^{(N-1)}).\label{eq-12}
\eeqa
Moreover, when $N\ge 2$ the functions  $F^i_j,G^i_j$ must satisfy the following equations:
\beq
F^j_{i,(k,N-1)}-F^j_{k,(i,N-1)}+(u^i-u^k)G^j_{k,(i,N-1)}=0. \label{eq-13}
\eeq
\end{lem}

\pf
By definition, we have
$$
Z^{ij}_{0,(k,2N+1)}=(-1)^{N+1}\left(f^i\,X^j_{(i,N)(k,N)}-g^i\,Y^j_{(i,N)(k,N)}\right).
$$
So vanishing of $Z^{ij}_{0,(k,2N+1)}$ implies that the functions $X^i, Y^i$ can be represented as
\begin{eqnarray*}
X^i&=&\sum_{j=1}^n(u^{j} {\tilde G}^i_{j}(u,\dots,u^{(N-1)};u^{j,N})+{\tilde F}^i_{j}(u,\dots,u^{(N-1)})u^{j,N})
+{\tilde Q}^i(u,\dots,u^{(N-1)}), \\
Y^i&=&\sum_{j=1}^n {\tilde G}^i_{j}(u,\dots,u^{(N-1)};u^{j,N}).
\end{eqnarray*}
When $N\ge 2$ the equations (\ref{eq-12}) of $X^i, Y^i$  follow from vanishing of
$$
Z^{ij}_{0,(i,2N)}=(-1)^{N+1}\left(N+\frac12\right)f^iu^i_xG^j_{i,(i,N)(i,N)}.
$$
In the case of $N=1$, the formulae (\ref{eq-12}) follow from vanishing of $Z^{ij}_{0,(i,2N)}$ and of
$Z^{ij}_2$.
Finally, for $N\ge 2$ the equation (\ref{eq-13}) is derived from vanishing of
$$
Z^{ij}_{0,(k,2N)}=(-1)^{N+1} f^i\left(F^j_{i,(k,N-1)}-F^j_{k,(i,N-1)}+(u^i-u^k)G^j_{k,(i,N-1)}\right).
$$
The lemma is proved.
\epf

\begin{lem}\label{lem-2}

(a) Assume that the vector fields $X, Y$ have components
of the form (\ref{eq-12}) and $N\ge 2$.
Then for any
$m=1,2,\dots,\left[\frac{N}2\right]$ the following identity holds true:
\beq\label{iden-1}
\sum_{l\ge0}(-1)^{m-l}\binom{N-l}{m-l}Z^{ij}_{l,(k,2N+1-m-l)}=(-1)^{N+1} f^i\,F^j_{i,(k,N-m)}.
\eeq
Here and below we denote $Z^{ij}$ as in (\ref{eq-z}) of Lemma \ref{lem-4}.\newline
\noindent (b) Assume that the vector field $X,Y$ have components
of the following form
\eqa
X^i&=&\sum_{j=1}^n \left(u^j G^i_j(u,\dots,u^{(N-m-1)};u^{j,N-m})+F^i_j(u,\dots,u^{(N-m-1)})\right)u^{j,N}\nn\\
&&\qquad+Q^i(u,\dots,u^{(N-1)}),\nn\\
Y^i&=&\sum_{j=1}^n G^i_j(u,\dots,u^{(N-m-1)};u^{j,N-m})\,u^{j,N}
+R^i(u,\dots,u^{(N-1)}) \label{eq-15}
\eeqa
with $N \ge 3$ and $1\le m\le \left[\frac{N-1}2\right]$. Then we have
\eqa
& & \sum_{l\ge0}(-1)^{m-l}\binom{N-l}{m-l}Z^{ij}_{l,(k,2N-m-l)} \nn \\
&=& (-1)^{N+1}\left[f^i\left(F^j_{i,(k,N-m-1)}-F^j_{k,(i,N-m-1)}\right)+f^i(u^i-u^k)G^j_{k,(i,N-m-1)}\right. \nn \\
& & \quad \left.
+(N-m+\frac12)f^iu^i_xG^j_{k,(i,N-m)}\dl^{ik}
+(A^{ik}-u^kB^{ik})G^j_{k,(k,N-m)}\right]\label{iden-2}
\eeqa
Here $A^{ik}, B^{ik}$ are defined by (\ref{eq-9}).
\end{lem}

\begin{prf}
By using the formula (\ref{lie-der}) the components $Z^{ij}$ of the bivector $\pal_1X-\pal_2Y$ have the explicit expressions
\eqa
Z^{ij}&=&\sum_{s\ge0}(-1)^{s+1}\left(f^i \p_x^{s+1}(X^j_{(i,s)}\dl)+\frac{\p_xf^i}2 \p_x^s(X^j_{(i,s)}\dl)
+\sum_{k=1}^n A^{ik} \p_x^s(X^j_{(k,s)}\dl)\right. \nn \\
&&\left.-g^i \p_x^{s+1}(Y^j_{(i,s)}\dl)-\frac{\p_xg^i}2 \p_x^s(Y^j_{(i,s)}\dl)-
\sum_{k=1}^n B^{ik} \p_x^s(Y^j_{(k,s)}\dl)\right)+\cdots. \nn
\eeqa
Here $\dl=\dl(x-y)$. It is easy to see that when $m\le\left[\frac{N}2\right]$, the first two terms in the formula
(\ref{lie-der}) don't appear in the identity (\ref{iden-1}). So we denote them by
periods
in the above formula and will omit them in the
calculations below. Then $Z^{ij}_p$ reads
\eqa
Z^{ij}_p&=&\sum_{s\ge0}(-1)^{s+1}\left(f^i\binom{s+1}{p}\p_x^{s+1-p}X^j_{(i,s)}+\frac{\p_xf^i}2
\binom{s}{p}\p_x^{s-p}X^j_{(i,s)}\right. \nn\\
&&\quad+\sum_{k=1}^n A^{ik}\binom{s}{p}\p_x^{s-p}X^j_{(k,s)}-g^i\binom{s+1}{p}\p_x^{s+1-p}\,Y^j_{(i,s)}\nn\\
&&\quad\left.-\frac{\p_xg^i}2 \binom{s}{p}\p_x^{s-p}\,Y^j_{(i,s)}
-\sum_{k=1}^n B^{ik}\binom{s}{p}\p_x^{s-p}\,Y^j_{(k,s)}\right)\nn
\eeqa
Denote by $l.h.s.$ the left hand side of the identity (\ref{iden-1}). We obtain
\eqa
&&l.h.s.= \nn\\
&&\sum_{p,s\ge0}(-1)^{m-p+s+1}\binom{N-p}{m-p}\left(
f^i\binom{s+1}{p}\sum_{t\ge0}\binom{s+1-p}{t}\p_x^tX^j_{(i,s)(k,2N-m-s+t)}\right.\nn\\
&&\qquad+\frac{\p_xf^i}2 \binom{s}{p}\sum_{t\ge0}\binom{s-p}{t}\p_x^tX^j_{(i,s)(k,2N+1-m-s+t)} \nn\\
&&\qquad+\sum_{l=1}^nA^{il}\binom{s}{p}\sum_{t\ge0}\binom{s-p}{t}\p_x^tX^j_{(l,s)(k,2N+1-m-s+t)}\nn\\
&&\qquad-g^i\binom{s+1}{p}\sum_{t\ge0}\binom{s+1-p}{t}\p_x^tY^j_{(i,s)(k,2N-m-s+t)}\nn\\
&&\qquad-\frac{\p_xg^i}2 \binom{s}{p}\sum_{t\ge0}\binom{s-p}{t}\p_x^tY^j_{(i,s)(k,2N+1-m-s+t)} \nn\\
&&\qquad\left.-\sum_{l=1}^nB^{il}\binom{s}{p}\sum_{t\ge0}\binom{s-p}{t}\p_x^tY^j_{(l,s)(k,2N+1-m-s+t)}\right). \nn
\eeqa
Here we used the commutation relations
$$
\frac{\p}{\p u^{i,q}}\p_x^m=\sum_{t\ge0}\binom{m}{t}\p_x^t\frac{\p}{\p u^{i,q-m+t}}
$$
By using the identity
$$
\sum_{p\ge0}(-1)^p\binom{N-p}{m-p}\binom{s}{p}\binom{s-p}{t}=\binom{s}{t}\binom{N-s+t}{m}
$$
and by changing the order of summation, we can rewrite $l.h.s.$ as follows
\eqa
&&l.h.s.= \nn\\
&&\sum_{s,t\ge0}(-1)^{m+s+1}\left(\binom{s+1}{t}\binom{N-s+t-1}{m}f^i\p_x^tX^j_{(i,s)(k,2N-m-s+t)}\right. \nn\\
&&\qquad+\binom{s}{t}\binom{N-s+t}{m}\frac{\p_xf^i}2 \p_x^tX^j_{(i,s)(k,2N+1-m-s+t)} \nn\\
&&\qquad+\binom{s}{t}\binom{N-s+t}{m}\sum_{l=1}^n A^{il}\p_x^tX^j_{(l,s)(k,2N+1-m-s+t)} \nn\\
&&\qquad-\binom{s+1}{t}\binom{N-s+t-1}{m}g^i\p_x^tY^j_{(i,s)(k,2N-m-s+t)} \nn\\
&&\qquad-\binom{s}{t}\binom{N-s+t}{m}\frac{\p_xg^i}2 \p_x^tY^j_{(i,s)(k,2N+1-m-s+t)} \nn\\
&&\qquad\left.-\binom{s}{t}\binom{N-s+t}{m}\sum_{l=1}^n B^{il}\p_x^tY^j_{(l,s)(k,2N+1-m-s+t)}\right)\label{form-1}
\eeqa
Now we substitute the expression (\ref{eq-12}) of $X,Y$ into the right hand side of the
above formula. By using the properties of binomial coefficients, it is easy to see that all terms in the above summation
vanish but the terms with $s = N, t = 0$, so the above formula can be simplified to
\eqa
l.h.s.&=&(-1)^{m+N+1}\left((-1)^m f^i X^j_{(i,N)(k,N-m)}-(-1)^m g^i Y^j_{(i,N)(k,N-m)}\right) \nn\\
&=&(-1)^{N+1} f^i F^j_{i,(k,N-m)} \nn
\eeqa
So part (a) of the lemma is proved.

Similarly to the derivation of (\ref{form-1}), we can prove the following formula:
\eqa
&&\sum_{l\ge0}(-1)^{m-l}\binom{N-l}{m-l}Z^{ij}_{l,(k,2N-m-l)}= \nn\\
&&\sum_{s,t\ge0}(-1)^{m+s+1}\left(\binom{s+1}{t}\binom{N-s+t-1}{m}f^i\p_x^tX^j_{(i,s)(k,2N-1-m-s+t)}\right. \nn\\
&&\qquad+\binom{s}{t}\binom{N-s+t}{m}\frac{\p_xf^i}2 \p_x^tX^j_{(i,s)(k,2N-m-s+t)} \nn\\
&&\qquad+\binom{s}{t}\binom{N-s+t}{m}\sum_{l=1}^n A^{il}\p_x^tX^j_{(l,s)(k,2N-m-s+t)} \nn\\
&&\qquad-\binom{s+1}{t}\binom{N-s+t-1}{m}g^i\p_x^tY^j_{(i,s)(k,2N-1-m-s+t)} \nn\\
&&\qquad-\binom{s}{t}\binom{N-s+t}{m}\frac{\p_xg^i}2 \p_x^tY^j_{(i,s)(k,2N-m-s+t)} \nn\\
&&\qquad\left.-\binom{s}{t}\binom{N-s+t}{m}
\sum_{l=1}^n B^{il}\p_x^tY^j_{(l,s)(k,2N-m-s+t)}\right).\label{form-2}
\eeqa
The summands of the r.h.s. vanish except for the terms with $(s,t)=(N,0),(N-m-1,0),(N-m,1),(N-m,0)$.
Then the $l.h.s$ of identity (\ref{iden-2}) reads
\eqa
l.h.s.&=&(-1)^{N+1}\left(f^iX^j_{(i,N)(k,N-m-1)}-g^iY^j_{(i,N)(k,N-m-1)}\right) \nn\\
&&+(-1)^N\left(f^iX^j_{(k,N)(i,N-m-1)}-g^iY^j_{(k,N)(i,N-m-1)}\right) \nn\\
&&+(-1)^{N+1}(N-m+1)\left(f^i\p_xX^j_{(k,N)(i,N-m)}-g^i\p_xY^j_{(k,N)(i,N-m)}\right) \nn\\
&&+(-1)^{N+1}\left(\frac{\p_xf^i}2X^j_{(k,N)(i,N-m)}-\frac{\p_xg^i}2Y^j_{(k,N)(i,N-m)}\right. \nn\\
&&\qquad\left.+\sum_{l=1}^n\left(A^{il}X^j_{(k,N)(l,N-m)}-B^{il}Y^j_{(k,N)(l,N-m)}\right)\right)\label{form-3}
\eeqa
Note that our $X,Y$ have properties
$$X^j_{(k,N)(i,N-m)}=u^kG^j_{k,(i,N-m)}\dl^{ki},\ Y^j_{(k,N)(i,N-m)}=G^j_{k,(i,N-m)}\dl^{ki}$$
Then the identity (\ref{iden-2}) follows from (\ref{form-3}) immediately. Part (b) of the lemma is proved.
\end{prf}

\begin{lem}\label{lem-ss-1}
Let $X, Y$ be two local vector fields that have components of the form (\ref{eq-11}) and satisfy the
relation
\beq\label{dxdy}
\pal_1 X=\pal_2 Y.
\eeq
Then the following statements hold true:\newline
\noindent i) When $N=2M+1$ the components of these vector fields have the expressions
\eqa
&&X^i=\sum_{j=1}^n {X}^i_j(u,u_x,\dots,u^{(M)}) u^{j,2M+1}+Q_i(u,u_x,\dots,u^{(2M)}),\label{ss-1}\\
&&Y^i=\sum_{j=1}^n {Y}^i_j(u,u_x,\dots,u^{(M)}) u^{j,2M+1}+R_i(u,u_x,\dots,u^{(2M)}).\label{ss-1-b}
\eeqa

\noindent ii) When $N=2M$ there exist local functionals $I_a, a=1,2,3$ such that the components of the
vector fields $X, Y$ have, after the modification
\beq\label{ss-5}
X \mapsto X-(\pal_1 I_1-\pal_2 I_2),\quad Y\mapsto Y-(\pal_1 I_2-\pal_2 I_3)
\eeq
(if necessary),
the expressions
\eqa
&&X^i=\sum_{j\ne i} {X}^i_j(u,u_x,\dots,u^{(M-1)}) u^{j,2M}+Q_i(u,u_x,\dots,u^{(2M-1)}),\label{ss-6}\\
&&Y^i=\sum_{j\ne i} {Y}^i_j(u,u_x,\dots,u^{(M-1)}) u^{j,2M}+R_i(u,u_x,\dots,u^{(2M-1)}).\label{ss-6-b}
\eeqa
In the case when the components $X^i$ of the vector field $X$ are homogeneous almost differential polynomials
of degree $d\ge 3$, we can choose the densities of the local functionals $I_a$ such that they
are homogeneous almost differential polynomials
of degree $d-1$.
\end{lem}

\pf For the case when $N=1$ the result of the lemma follows from Lemma \ref{lem-1},
 so we assume that $N\ge 2$.
It follows from Lemma \ref{lem-1} that the components of $X,Y$
must take the form (\ref{eq-12}). The result of part (a) of Lemma \ref{lem-2} then shows that $F^i_j$ are
independent of $u^{(N-m)}$ for $m=1,\dots,\left[\frac{N}2\right]$, so the identities in (\ref{eq-13}) read
$$
(u^i-u^k)G^j_{k,(i,N-1)}=0,
$$
thus $G^j_k$ are independent of $u^{i,N-1}$ when $i \ne k$. When $N\ge 3$ we use the identity (\ref{iden-2})  of
Lemma \ref{lem-2} with $m=1$ to obtain, by putting $i=k$,
$$
G^j_{i,(i,N-1)}=0.
$$
So $G^j_k$ are in fact independent of $u^{(N-1)}$.
Now the identity (\ref{iden-2}) shows that
$$
G^j_{k,(i,N-2)}=0,\quad i\ne k.
$$
By using repeatedly the identity (\ref{iden-2}) we know that $G^i_j$ is independent of $u^{(N-m)}$
for $m=1,\dots,\left[\frac{N-1}2\right]$, and $G^j_{i,(k,M)}=0$ for $N=2M\ge 2$ and $i\ne k$.

For the case of $N=2M+1, M\ge 1$ the above argument shows that the components of the vector
fields $X,Y$ have the form (\ref{ss-1}), (\ref{ss-1-b})
where ${X}^i_j, {Y}^i_j$ are given by the expressions
\eqa
&&{X}^i_j=u^j G^i_j(u,u_x,\dots,u^{(M)})+F^i_j(u,u_x,\dots,u^{(M)}),\nn\\
&&{Y}^i_j= G^i_j(u,u_x,\dots,u^{(M)}.
\eeqa

For the case of $N=2M, M\ge 1$ the above argument shows that $X^i, Y^i$
have the form
\eqa
X^i&=&\sum_{j=1}^n \left(u^j G^i_j(u,\dots,u^{(M-1)},u^{j,M})+F^i_j(u,\dots,u^{(M-1)})\right)u^{j,N}\nn\\
&&\qquad+Q^i(u,\dots,u^{(N-1)}),\label{eq-15b} \\
Y^i&=&\sum_{j=1}^n G^i_j(u,\dots,u^{(M-1)},u^{j,M})\,u^{j,N}
+R^i(u,\dots,u^{(N-1)}) \label{eq-15bb}
\eeqa
{}From vanishing of the coefficients of $\delta^{(2M+1)}(x-y)$ in the expression of $(\pal_1 X-\pal_2 Y)^{ij}$
it follows that
\beq\label{ss-3}
f^i(u) F^j_i+f^j(u) F^i_j=0.
\eeq
In particular, we have
\beq\label{ss-4}
F^i_i=0, \quad i=1,\dots,n.
\eeq
Define the functionals
\beq\label{ss-2}
I_k=(-1)^M \int\sum_{i=1}^n \pal_{u^{i,M}}^{-2} \frac{(u^i)^{3-k}\,G^i_i}{(M+\frac12) f^i u^i_x} dx,\quad k=1,2,3
\eeq
Then the components of the vector fields $\tilde X=X-(\pal_1 I_1-\pal_2 I_2), \ \tilde Y= Y-(\pal_1 I_2-\pal_2 I_3)$
have the
form of (\ref{eq-15b}), (\ref{eq-15bb}) with $G^i_i=0$. Since the vector fields $\tilde X, \tilde Y$ still satisfy the
relation $\pal_1 \tilde X=\pal_2 \tilde Y$, we can assume without
loss of generality that the components of
$X, Y$ have the form (\ref{eq-15b}), (\ref{eq-15bb}) with vanishing $G^i_i$.
By using the equations (\ref{ss-3}) we check
that the identity
(\ref{iden-2}) is still valid for $i=k$ when $m=M$ and $N\ge 2$.
This leads to the fact that
the functions $G^i_j$ for $i\ne j$ do not depend on $u^{j,M}$. Thus
we proved that
the components of $X, Y$ have the form (\ref{ss-6}), (\ref{ss-6-b})
after the modification (\ref{ss-5}) if necessary.

When the components of the vector field $X$ are homogeneous almost differential polynomials of degree
$d\ge 3$, the equations (\ref{ss-4}) and the expression (\ref{eq-15b}) imply that the fuctions $G^i_i$ are
also  homogeneous almost differential polynomials of degree
$d-2 M$. So, when $M\ge 2$ we can choose the densities of the functionals $I_a$ defined in (\ref{ss-2}) to be
homogeneous almost differential polynomials of degree $d-1$. In the case $M=1$, since the functions
$G^i_i=G^i_i(u;u^i_x)$ are homogeneous of degree $d-2\ge 1$ (recall our assumption $d\ge 3$),
the function $\frac{G^i_i(u;u^i_x)}{u^i_x}$
is in fact a polynomial in $u^i_x$. Thus in this case we can still choose the densities of the
functionals $I_a$ defined in (\ref{ss-2}) to be
homogeneous almost differential polynomials of degree $d-1$. The lemma is proved.\epf

\begin{lem}\label{lem-ss-2}
Let the vector fields $X, Y$ have components of the form (\ref{eq-11}) with $N\ge 2$
and satisfy the relation (\ref{dxdy}).
If the functions $X^i, i=1,\dots,n$ do not depend on $u^{(N)}$, then we can modify the vector field $Y$
by
\beq
Y\mapsto Y-\pal_2 J
\eeq
for certain local functional $J$ such that the components of this modified vector field $Y$
depend at most on $u,\dots, u^{(N-1)}$ and the relation (\ref{dxdy}) still holds true.
\end{lem}
\pf
We first assume that
$N=2M+1$. From the assumption of the lemma and the result of Lemma \ref{lem-ss-1} we know that
the components of the vector fields $X, Y$ have the form (\ref{ss-1}), (\ref{ss-1-b}) with $X^i_j=0$.
To prove the lemma we need to find a local functional $J$ with density $h(u,u_x,\dots,u^{(M)})$
which satisfies the conditions
\beq\label{ss-7}
\frac{\pal^2 h}{\pal u^{i,M}\pal u^{j,M}}=(-1)^M \frac1{g^i(u)} Y^i_j,\quad i, j=1,\dots,n.
\eeq
Denote by $A_{ij}$ the r.h.s. of the above formulae. Then
from vanishing of the coefficients of
$\delta^{(2M+2)}(x-y)$ in the expression of the components of the bivector $\pal_1 X-\pal_2 Y$ it follows that
the functions $A_{ij}$ are symmetric with respect to the indices $i,j$.
{}From the equation (\ref{iden-2}) with $m=M$ we also know that
\beq
\frac{\pal Y^j_i}{\pal u^{k,M}}=\frac{\pal Y^j_k}{\pal u^{i,M}}.
\eeq
So the functions $\frac{\pal A_{ij}}{\pal u^{k,M}}$ are symmetric with respect to the indices $i,j,k$ which
implies the existence of a function $h(u,u_x,\dots,u^{(M)})$ satisfying the requirement (\ref{ss-7}).

Next let us assume that $N=2M$. As we did in the proof of Lemma \ref{lem-ss-1} we can show that the components
of $X, Y$ have the form (\ref{eq-15b}), (\ref{eq-15bb}) with
\beq\label{ss-10}
F^i_j=-u^j G^i_j.
\eeq
Since the functions $F^i_j$ do not depend on $u^{(M)}$
we deduce that the functions $X^i, Y^i$
must have the expressions (\ref{ss-6}) with $X^i_j=0$.
{}From (\ref{ss-4}) and the independence of $X^i$ on
$u^{(N)}$ it also follows  that $G^i_i=0$ for $i=1,\dots,n$.  Now by using the
vanishing of the coefficients of $\delta^{(2M+1)}$ of the components of $\pal_1X-\pal_2 Y$ and that
of the l.h.s. of (\ref{iden-2})
with $m=M$ we obtain
\eqa\label{eq-25}
&&\hy_{ij}+\hy_{ji}=0,\\
&&\hy_{jk,(i,M-1)}-\hy_{ik,(j,M-1)}-
\hy_{ji,(k,M-1)}=0.\label{eq-24}
\eeqa
Here $\hy_{ij}=\frac1{g^i} Y^i_j$. The above two equations ensure the existence of a 1-form
$\al=\sum_{i=1}^n h_i(u,\dots,u^{(M-1)}) d u^{i,M-1}$ such that
\beq
d\al=\frac12 \sum_{i,j} \hy_{ij}\, d u^{i,M-1}\wedge d u^{j,M-1}.
\eeq
Now the functional $J$ defined by
\beq
J=\int \sum_{i=1}^n h_i(u(x),\dots,u^{(M-1)}(x))\, u^{i,M}(x) dx
\eeq
meets the requirement of the lemma and we finished the proof.\epf

\begin{lem}\label{lem-3}
Let $X$ be a local vector field with components
$$
X^i=X^i(u;u_x,\dots,u^{(N)}),\quad i=1,\dots,n,\ N\ge 4.
$$
that are homogeneous almost differential polynomials
of degree $d\ge 3$. If $X$ also has the following properties:\newline
\noindent (a) When $N=2M+2$, the  components of $X$ have the form
\beq
X^i=\sum_{j\ne i} X^i_j(u,\dots,u^{(M)})u^{j,2M+2}+Q^i(u,\dots,u^{(2M+1)})
\label{eq-18}
\eeq
and satisfy the conditions
\beq
(u^k-u^j)\left(\hx_{ij,(k,M)}-\hx_{ik,(j,M)}\right)+(u^k-u^i)\hx_{jk,(i,M)}+
(u^j-u^i)\hx_{kj,(i,M)}=0
\label{cond-3}
\eeq
for
\beq\label{ss-11}
\hx_{ij}=\frac1{f^i(u)}{X^i_j(u,\dots,u^{(M)})}.
\eeq
\noindent (b) When $N=2M+1$, $X$ has components of the form
\beq
X^i=\sum_{j=1}^n X^i_j(u,\dots,u^{(M)})u^{j,2M+1}+Q^i(u,\dots,u^{(2M)}) \label{eq-16}
\eeq
and satisfies the conditions
\eqa
&&X^i_{j,(k,M)}-X^i_{k,(j,M)}=0 \label{cond-1}\\
&&(u^k-u^j)\hx_{ij,(k,M)}+(u^i-u^k)\hx_{jk,(i,M)}+(u^j-u^i)\hx_{ki,(j,M)}=0.
\label{cond-2}
\eeqa
Then there exist two local functionals $I_1,I_2$ with densities that are homogeneous
almost differential polynomials of degree $d-1$
such that the components of the vector field $X-(\pal_1 I_1-\pal_2 I_2)$ depend at most on
$u,u_x,\dots,u^{(N-1)}$.
\end{lem}

\begin{prf}
We first prove the lemma for the case when $d\ge 3, N\ge 5$.
Assume $N=2M+2$ and the vector field $X$ satisfies the conditions (\ref{eq-18}), (\ref{cond-3}). We want to
find two local functionals $I_1, I_2$ with densities of the form
\beq\label{eq-27}
h_a=\sum_{j=1}^n h_{a;j}(u,u_x,\dots,u^{(M)}) u^{j,M+1},\quad a=1,2
\eeq
such that they meet the requirement of the lemma. For this we need
to find the functions
$h_{a;j}, a=1,2, j=1,\dots, n$
satisfying the following equations:
\beq
(-1)^{M+1} X^i_j=f^i\left(\frac{\pal h_{1;i}}{\pal u^{j,M}}-\frac{\pal h_{1;j}}{\pal u^{i,M}}\right)-
g^i\left(\frac{\pal h_{2;i}}{\pal u^{j,M}}-\frac{\pal h_{2;j}}{\pal u^{i,M}}\right).\label{eq-23}
\eeq

Denote
\beq
P_{ij}=(-1)^{M} \frac{g^j X^i_j+g^i X^j_i}{f^i g^j-f^j g^i},\quad
Q_{ij}=(-1)^{M} \frac{f^j X^i_j+f^i X^j_i}{f^i g^j-f^j g^i}
\eeq
Then it follows from (\ref{cond-3}) that the two-forms
\beq\label{2-form}
\varpi_1=\frac12 \sum_{i,j} P_{ij} du^{i,M}\wedge du^{j,M},\quad
\varpi_2=\frac12 \sum_{i,j}Q_{ij} du^{i,M}\wedge du^{j,M}
\eeq
are closed. So there exist one-forms
\beq\label{1-form}
\al_a=\sum_{j} h_{a;j}(u,u_x,\dots,u^{(M)})\, du^{j,M},\ a=1,2
\eeq
such that $d \al_a=\varpi_a$ and the functions $h_{a;j}$ are homogeneous
almost differential polynomials of degree $d-M-2$. Now it's easy to see that the functions $h_{a;j}$ satisfy
(\ref{eq-23}). So we proved the lemma for $N=2M+2>4$.

Next we assume that $N=2M+1>4$ and the vector field $X$ satisfy the condition (\ref{eq-16})-(\ref{cond-2}). Let the
two local functionals $I_1, I_2$ that we are looking for have densities of the form
\beq\label{eq-26}
h_a(u,u_x,\dots,u^{(M)}),\quad  a=1,2.
\eeq
Since
the $i$-th component of the vector field $\pal_1 I_1-\pal_2 I_2$ depends
at most on $u$,$\dots$,$u^{(2M+1)}$ and, moreover,
it depends linearly on $u^{k,2M+1}$, we only need to find functions $h_1, h_2$ such that
\beq\label{eq-19}
X^i_j=\frac{\pal (\pal_1 I_1-\pal_2 I_2)^i}{\pal u^{j,2M+1}}=(-1)^M (f^i\, h_{1,(i,M)(j,M)}-g^i\, h_{2,(i,M)(j,M)}).
\eeq
To this end, let us define $P_{ii}, Q_{ii},\ i=1,\dots,n$
by solving the following systems
\eqa
&&\frac{\pal P_{ii}}{\pal u^{j,M}}=
\frac{(-1)^{M+1}}{u^i-u^j}\left(u^j \frac{X^i_{j,(i,M)}}{f^i}-u^i \frac{X^j_{i,(i,M)}}{f^j}\right),\quad j\ne i,\nn\\
&&\frac{\pal Q_{ii}}{\pal u^{j,M}}=
\frac{(-1)^{M+1}}{u^i-u^j}\left(\frac{X^i_{j,(i,M)}}{f^i}-\frac{X^j_{i,(i,M)}}{f^j}\right),\ j\ne i.\label{eq-21}
\eeqa
The conditions (\ref{cond-1}), (\ref{cond-2}) implies the compatibility of the above systems, i.e.,
\beq
\frac{\pal}{\pal u^{k,M}}\left(\frac{\pal P_{ii}}{\pal u^{j,M}}\right)
=\frac{\pal}{\pal u^{j,M}}\left(\frac{\pal P_{ii}}{\pal u^{k,M}}\right),\quad j,k\ne i.
\eeq
So we have a set of functions $P_{ii}=P_{ii}(u,u_x,\dots,u^{(M)}),
Q_{ii}=Q_{ii}(u,u_x,\dots,u^{(M)})$ satisfying the conditions (\ref{eq-21}). The ambiguity in
the definition of these functions is the following shifts:
\eqa
&&P_{ii}\mapsto P_{ii}+W_{ii}(u,u_x,\dots,u^{(M-1)},u^{i,M}),\nn\\
&&Q_{ii}\mapsto Q_{ii}+R_{ii}(u,u_x,\dots,u^{(M-1)},u^{i,M}).\label{eq-22}
\eeqa
Here $W_{ii}, R_{ii}$ are arbitrary functions to be specified later. We also define functions
$P_{ij}, Q_{ij}$ with $i\ne j$ by the following formulae:
\beq
P_{ij}=\frac{(-1)^{M+1}}{u^i-u^j}\left(u^j \frac{X^i_j}{f^i}-u^i \frac{X^j_i}{f^j}\right),
\ Q_{ij}=\frac{(-1)^{M+1}}{u^i-u^j}\left(\frac{X^i_j}{f^i}-\frac{X^j_i}{f^j}\right).
\eeq
By using the conditions (\ref{cond-1}), (\ref{cond-2}) we easily verify that
$$
\frac{\pal P_{ij}}{\pal u^{k,M}},\quad \frac{\pal Q_{ij}}{\pal u^{k,M}}
$$
are symmetric with respect to their indices $i,j,k$.
So there exist functions \newline
$h_a(u,u_x,\dots,u^{(M)}),\ a=1,2$ such that
\beq
\frac{\pal^2 h_1}{\pal u^{i,M}\pal u^{j,M}}=P_{ij},\quad \frac{\pal^2 h_2}{\pal u^{i,M}\pal u^{j,M}}=Q_{ij}.
\eeq
Now it's easy to verify that, when $i\ne j$, the functions $h_1, h_2$
satisfy the conditions given in (\ref{eq-19}).
When $i=j$ we have
\beq
X^i_i-(-1)^M (f^i\, h_{1,(i,M)(i,M)}-g^i\, h_{2,(i,M)(i,M)})
=X^i_i-(-1)^M (f^i\, P_{ii}-g^i\, Q_{ii}).
\eeq
It follows from the defintion of $P_{ii}, Q_{ii}$ and the condtions given in (\ref{cond-1}) that
the r.h.s. of the above formulae does not depend on $u^{k,M}$ for any $k\ne i$, so we can make it to
be zero by adjusting the functions $P_{ii}, Q_{ii}$ as in (\ref{eq-22}).
By the above construction,
the functions $h_1, h_2$ can be chosen to be homogeneous
almost differential polynomials of degree $d-1$.
So  the lemma is proved for the case mentioned above.

Now let us consider the case $d\ge 3, N=4$. Proceeding in the same way as for the case of $N=2M+2, M\ge 2$ we
can find the 1-forms (\ref{1-form}) such that the 2-forms $\varpi_1, \varpi_2$ that are defined as in
(\ref{2-form}) can be represented by $\varpi_a=d \al_a,\ a=1,2$. The pecularity of this particular case of $M=1$
lies in the fact that the functions $h_{a,j}$ that we constructed above are in general no longer rational functions
of the jet coordinates $u^{i,k}, k\ge 1$, they can be chosen to have the form
\beq\label{h-aj}
h_{a,j}=\sum_{k=1}^n W_{a,j;k}(u;u_x) \log u^k_x+U_{a,j}(u;u_x),\quad a=1,2,\ j=1,\dots,n.
\eeq
Here $W_{a,j;k}, U_{a,j}\in {\cal A}$ and are homogeneous of degree $d-1$. Since
\beq
\frac{\pal h_{a,j}}{\pal u^i_x}-\frac{\pal h_{a,i}}{\pal u^j_x}\in {\cal A}
\eeq
we must have
\beq
\frac{\pal W_{a,j;k}}{\pal u^i_x}=\frac{\pal W_{a,i;k}}{\pal u^j_x},\quad i, j=1,\dots, n.
\eeq
This implies the existence of  functions $A_{a,k}(u;u_x)\in {\tilde A}$ of degree
$d-2$ such that
\beq
W_{a,j;k}=\frac{\pal A_{a,k}(u;u_x)}{\pal u^j_x},\quad a=1,2;\ j,k=1,\dots,n.
\eeq
Since $W_{a,j;k}$ are almost differential polynomials, the functions $A_{a,k}$ can also be chosen as homogeneous
almost differential polynomials of degree $d-2\ge 1$
up to the addition of terms of the form
\beq
\sum_{l=1}^n B_{a,k,l}(u) \log u^l_x.
\eeq
However, such functions have degree zero, so they are not allowed to appear in the expression of $A_{a,k}$.
Now the needed functionals $I, J$ that satisfy the requirements of the lemma can be chosen to have
densities ${\tilde h}_1, {\tilde h}_2$ of the form
\beq
\tilde h_a=\sum_{k=1}^n\left(U_{a,k}(u,u_x)- A_{a,k}(u,u_x) \frac{1}{u^k_x}\right) u^k_{xx},\quad a=1,2.
\eeq
The lemma is proved.
\end{prf}
\vskip 0.3cm

\noindent{\em Proof of Theorem \ref{thm3}.}
Assume that the components $X^i$ of the vector field $X$ have the form
\beq\label{ss-8}
X^i=X^i(u,\dots,u^{(N)}),\quad N\ge 4,\ i=1,\dots,n.
\eeq
The relation (\ref{d1d2}) implies the existence of a local vector field $Y$ such that $\pal_1 X=\pal_2 Y$.
By using Lemma \ref{lem-ss-2} we can choose the vector field $Y$ such that
its components depend
at most on the coordinates $u,\dots,u^{(N)}$. Then it follows from the results of Lemma \ref{lem-ss-1}
that the components of $X, Y$ have the expressions (\ref{ss-1}),(\ref{ss-1-b}) when $N$ is odd and the one of
(\ref{ss-6}),(\ref{ss-6-b})
when $N$ is even (after a modification of (\ref{ss-5}) which does not affect our
result). We now proceed to employ the result of Lemma \ref{lem-3} in order to
find two local functionals $I, J$ with densities that are
almost differential polynomials of degree $d-1$ such that $X-(\pal_1I -\pal_2 J)$ depends at most on $u,\dots,u^{(N-1)}$.
To this end we need to verify that $X^i$ satisfy the equations (\ref{cond-3}) when $N$ is even and
that of (\ref{cond-1}), (\ref{cond-2}) when $N$ is odd.

Let $N$ be an even integer $2M+2$. Then by using vanishing of the l.h.s. of (\ref{iden-2})
with $m=M$  we obtain
\beq\label{eq-24-b}
{\hat X}_{jk,(i,M)}-\hx_{ik,(j,M)}-\hx_{ji,(k,M)}=u^i u^j \left(\hy_{jk,(i,M)}- \hy_{ik,(j,M)}-
\hy_{ji,(k,M)}\right).
\eeq
Here $\hy_{ij}$ are defined as in (\ref{eq-24}) and $\hx_{ij}$ are defined by (\ref{ss-11}).
We also have the following equations
\beq\label{eq-25-b}
\hx^j_i+\hx^i_j=u^i u^j \left( \hy_{ji}+\hy_{ij}\right)
\eeq
due to vanishing of the coefficients of $\delta^{(2M+1)}$ in the components of $\pal_1X-\pal_2 Y$.
Denote by $L_{i,j,k}, R_{i,j,k}$ the expression of the l.h.s. and r.h.s. of (\ref{eq-24-b}) multiplied by
$u^k$. Then we have
\eqa
&&R_{i,j,k}-R_{j,k,i}-R_{k,i,j}-R_{i,k,j}=u^i u^j u^k \left(\hy_{jk,(i,M)}+\hy_{kj,(i,M)}\right)\nn\\
&&=u^i \left( \hx_{jk,(i,M)}+\hx_{kj,(i,M)}\right).\nn
\eeqa
Here in the last equality we used the equations in (\ref{eq-25-b}). By equating the last
expression with $L_{i,j,k}-L_{j,k,i}-L_{k,i,j}-L_{i,k,j}$ we arrive at the proof that the vector field $X$
satisfies the conditions given in (\ref{cond-3}).

Now let us assume that
$N=2M+1$. Then from (\ref{iden-2}) with $m=M$  we know that
\beq\label{ss-12}
f^i(X^j_{i,(k,M)}-X^j_{k,(i,M)})=g^i(Y^j_{i,(k,M)}-Y^j_{k,(i,M)}),
\eeq
This equation together with the one that is obtained from it by exchange the indices $i$ and $k$ implies that
\beq
X^j_{i,(k,M)}=X^j_{k,(i,M)},\quad Y^j_{i,(k,M)}=Y^j_{k,(i,M)}.
\eeq
So the components $X^i, Y^i$ satisfy the conditions of the form (\ref{cond-1}). We are left to prove that
they also satisfy the conditions of the form (\ref{cond-2}). For this,
let us consider the coefficients
of $\delta^{(2M+2)}(x-y)$  in the expression of the components
of the bivector $\pal_1 X-\pal_2 Y$.
Vanishing of these coefficients leads to the equations
\beq\label{ss-13}
{\hat X}_{ji}-{\hat X}_{ij}=u^i u^j \left({\hat Y}_{ji}-{\hat Y}_{ij}\right).
\eeq
By taking derivative with respect to $u^{k,M}$ and multiplying by $u^k$ on both sides of the above equation we
obtain
\beq
u^k\left({\hat X}_{ji,(k,M)}-{\hat X}_{ij,(k,M)}\right)=u^k u^i u^j \left({\hat Y}_{ji,(k,M)}-
{\hat Y}_{ij,(k,M)}\right).
\eeq
Denote the r.h.s. of the last equation by $W_{i,j,k}$. Then the condition (\ref{cond-2}) follows from
$W_{i,j,k}+W_{j,k,i}+W_{k,i,j}=0$.

Above we showed that the vector field $X$ satisfies the requirements of
Lemma \ref{lem-3}. So we can find
local functionals $I, J$ with densities that are
almost differential polynomials of degree $d-1$ such that $X-(\pal_1I -\pal_2 J)$ depends at most on $u,\dots,u^{(N-1)}$.
Repeating this procedure by the subtraction of terms of the form $\pal_1 I-\pal_2 J$, we reduce the proof of the Theorem
\ref{thm3} to the case when the components $X^i$ of the vector field $X$ have the form (\ref{ss-8}) with $N=1,2,3$.

Note that in the case when $N=3$ or $N=2$, the components of the vector field $X$ and the accompanying one $Y$
also have the forms (\ref{ss-1}),(\ref{ss-1-b}) and
(\ref{ss-6}),(\ref{ss-6-b}) with $M=1$, and the above equations (\ref{eq-24-b}), (\ref{eq-25-b}) and
(\ref{ss-12}), (\ref{ss-13}) still hold true. Thus when $N=3$
the vector field $X$ fulfils the requirements of Lemma \ref{lem-3},
and we can find local functionals $I_1, I_2$ such that $X-(\pal_1 I_1-\pal_2 I_2)$
depend at most on $u,u_x,u_{xx}$. The difference of this special case from
the general one
lies in the fact that now it is not obvious that we can choose the densities $h_a(u,u_x), a=1,2$ to be
almost differential polynomials. What can easily be seen from
our construction is that they can be chosen to have the form
\beq\label{ss-9}
h_a(u,u_x)=\sum_{i\ne j} V_{a;i,j}(u,u_x) \log u^i_x\log u^j_x+\sum_{i=1}^n V_{a;i}(u,u_x) \log u^i_x+U_a(u,u_x).
\eeq
Here the functions $V_{a;i,j}, V_{a;i}, U_a$ are homogeneous almost differential polynomials of degree $d-1$.
The vector field $\tilde X=X-(\pal_1 I_1-\pal_2 I_2)$ still has the property $\pal_1 \pal_2 \tilde X=0$,
so the same argument as above shows that
the components of $\tilde X$ have the form (\ref{eq-18}) with $M=0$ and satisfy
the equations (\ref{cond-3}). By using the construction of Lemma \ref{lem-3} we can find local functionals $I_3, I_4$
such that the vector field $\tilde X-(\pal_1 I_3-\pal_2 I_4)$ depends at most on $u, u_x$. A careful analysis of this
construction shows that the densities of these local functionals can also be chosen to have the form of
(\ref{ss-9}). Let us denote the densities $I=I_1+I_3,\ J=I_2+I_4$ also by $h_1, h_2$ which have the expression
(\ref{ss-9}). A simple calculation shows that
\eqa
&&\frac{\pal(\pal_1 I-\pal_2 J)^i}{\pal u^{j}_{xxx}}=g^i\frac{\pal^2 h_2}{\pal u^i_x\pal u^j_x}-
f^i\frac{\pal^2 h_1}{\pal u^i_x\pal u^j_x},\\
&&
\frac{\pal (\pal_1 I-\pal_2 J)^i}{\pal u^i_{xx}}=-\frac32 f^i(u) u^i_x \frac{\pal^2 h_2}{\pal u^i_x \pal u^i_x}.
\eeqa
The r.h.s. of the above two identities equal respectively to $\frac{\pal X^i}{\pal u^j_{xxx}}$ and
$\frac{\pal X^i}{\pal u^j_{xx}}$. We deduce that the functions
\beq
\frac{\pal^2 h_a}{\pal u^i_x \pal u^j_x},\quad a=1,2;\ i,j=1,\dots,n
\eeq
are homogeneous almost differential polynomials of degree $d-3$. This fact yields the restriction on the
coefficients $V_{a;i,j}, V_{a;i}$ in the expression of the densities (\ref{ss-9}) that they depend on $u_x$
at most linearly. Since $h_1, h_2$ have degree $d-1\ge 2$ (recall that we assume $d\ge 3$),
it follows that the functions $V_{a;i,j}, V_{a;i}$
must vanish, and as a result the densities $h_1, h_2$ of the local functionals $I, J$ are homogeneous almost differential polynomials of degree $d-1$. Now let us prove that we have in fact
\beq
X=\pal_1 I-\pal_2 J.
\eeq
This is due to the fact that the vector field $\bar X=X-(\pal_1 I-\pal_2 J)$ still satisfies
the property $\pal_1 \pal_2
\bar X=0$. So, by using Lemma \ref{lem-ss-2} we can find a vector field $\bar Y$ that depends at most on
$u,u_x$ such that $\pal_1 \bar X=\pal_2 \bar Y$. Then by using Lemma \ref{lem-1} we know that the components
of $\bar X$ depend at most linearly on $u_x$, since they are homogeneous almost differential polynomials
of degree $d\ge 3$ we must have $\bar X=0$. Thus we proved the theorem.\epf

\begin{thm}\label{tmd-d2}
Let the vector field $X$ satisfying the condition (\ref{d1d2}) have components of the form
\beq\label{ss-14}
X^i=\sum_{j=1}^n X^i_j(u) u^j_{xx}+\sum_{k,l}Q^i_{kl}(u) u^k_x u^l_x.
\eeq
Then each function $\hx_{ii}=(f^i(u))^{-1} X^i_i(u)$ depends only on $u^i$ and
there exist local functionals $\tilde I_1, \tilde I_2$
with densities that are homogeneous
differential polynomials of degree $1$ such that
\beq
X=\pal_1 (I_1+\tilde I_1)-\pal_2 (I_2+\tilde I_2)
\eeq
where
\beq\label{ss-15}
I_a=-\frac23\sum_{i=1}^n \int (u^i)^{2-a} \hx_{ii}(u^i) u^i_x \log u^i_x dx
,\quad a=1,2.
\eeq
\end{thm}

Note that the densities $h_1, h_2$ of the functionals $I_1, I_2$ can be chosen as
\beq
h_a=\sum_{i=1}^n V_{a,i}(u^i) \frac{u^i_{xx}}{u^i_x}\, dx,\quad a=1,2
\eeq
which are homogeneous almost differential polynomials of degree $d-1=1$.
Here the functions $V_{a,i}$ are defined by $V_{a,i}(u^i)'=2 (u^i)^{2-a} \hx_{ii}(u^i)$.
\vskip 0.2cm
\begin{prf}
Since $X^i$ are differential polynomials the condition (\ref{d1d2}) implies existence of a vector field $Y$ with components $Y^i$ of the same form (\ref{ss-14}) of $X^i$ such that
$\pal_1 X=\pal_2 Y$. By using vanishing of
the l.h.s. of (\ref{iden-2})  with $k=i\ne j$ we deduce that
\beq\label{ff}
\frac{\pal \hx_{ii}}{\pal u^j}-u^iu^j\frac{\pal \hy_{ii}}{\pal u^j}=0,\quad i\ne j=1,\dots,n.
\eeq
Since (\ref{eq-25-b}) also holds true in this case, we obtain
\beq
\hy_{ii}=\frac1{(u^i)^2} \hx_{ii}
\eeq
which yields, together with (\ref{ff}), the first result of the theorem
\beq
\frac{\pal \hx_{ii}}{\pal u^j}=0,\quad i\ne j=1,\dots,n.
\eeq

{}From the definition (\ref{ss-1}) it is easy to see that the components of the vector field
$\tilde X=X-(d_1 I_1-d_2 I_2)$ are still homogeneous differential polynomials
of the form
(\ref{ss-4}) with $X^i_i(u)=0,\ i=1,\dots, n$. Then by using the same construction as we give
in the proof of Lemma \ref{lem-3}  for the local functionals with the densities (\ref{eq-27})
we can find functionals $\tilde I_1, \tilde I_2$ with homegeneous differential
polynomial densities of degree $1$
such that the vector field $\bar X=\tilde X-(\pal_1 \tilde I_1
-\pal_2 \tilde I_2)$ depends at most on $u, u_x$.
The equation $\pal_1 \pal_2 \bar X=0$ then implies that $\bar X$ depends
at most linearly on $u_x$.
Since the components of $\bar X$ are homogeneous differential polynomials
of degree 2, we
arrive at the equalities
\beq
X=\pal_1 (I_1+\tilde I_1)-\pal_2 (I_2+\tilde I_2).
\eeq
The Theorem is proved.
\end{prf}

\vskip 0.2cm
\noindent{\em Proof of the Quasitriviality Theorem}\quad
Due to the triviality of the Poisson cohomology\newline
$H^2({\cal L}(M),\varpi_1)$ the bihamiltonian
structure (\ref{eq-2}) can always be assumed, if necessary by performing a usual Miura-type transformation,
to have the  following form:
\begin{eqnarray}
\{u^i(x),u^j(y)\}_1&=&\{u^i(x),u^j(y)\}^{[0]},\\
\{u^i(x),u^j(y)\}_2&=&\{u^i(x),u^j(y)\}^{[0]}
+\sum_{k \ge 1}\ve^k Q^{ij}_k.\label{eq-20b}
\end{eqnarray}
Here $Q^{ij}_k$ are the components of the bivectors $Q_k$ and have the expressions
$$
Q^{ij}_k=\sum_{l=0}^{k+1}Q^{ij}_{k,l}(u;u_x,\dots,u^{(k+1-l)})\,\delta^{(l)}(x-y).
$$
We also denote by $Q_0$ the bivector corresponding to the undeformed second Poisson structure.
The coefficients $Q^{ij}_{k,l}$ are homogeneous differential polynomials of degree $k+1-l$.
The compatibility of the above two Poisson brackets implies the existence of vector fields $X_k,\ k\ge 1$
such that
\beq
Q_k=\pal_1 X_k, \quad \pal_1 \pal_2 X_1=0,\quad k\ge 0
\eeq
and the components of $X_k$ are homogeneous differential polynomials of degree $k$.
By using Theorem \ref{thm3}, we know the existence of two local functionals
\beq
I_1=\int h_{1,1}(u(x)) dx,\quad J_1=\int h_{2,1}(u(x)) dx
\eeq
such that $X_1=\pal_1 I_1-\pal_2 J_1$. By performing the Miura-type transformation
\beq
u^i\mapsto \exp({\ve\,{\pal_1 J_1}})\, u^i
\eeq
the first Poisson structure remains the same while the second Poisson structure (\ref{eq-20b}) is
transformed to
\beq\label{ss-20}
\exp({-\ve\, \ad_{\pal_1 J_1}})(Q_0+\sum_{k\ge 1} \ve^k \pal_1 {X}_k)=
Q_0+\sum_{k\ge 2} \ve^k \pal_1 {\tilde X}_k
\eeq
where the vector fields $\tilde X_k$ have the
expressions
\beq
{\tilde X}_k=\sum_{l=0}^k (-1)^l \frac{\left(\ad_{\pal_1 J_1}\right)^l}{l!} X_{k-l},\quad k\ge 2.
\eeq
The components of these vector fields are homogeneous differential polynomials of degree $k$.
So we can use Theorem \ref{thm3} again to find two local functionals $I_2, J_2$ with densities $h_{a,2}(u,u_x,u_{xx})$
that are homogeneous almost differential polynomials of degree 1 such that $X_2=\pal_1 I_2-\pal_2 J_2$. Then the Miura-type transformation
\beq
u\mapsto \exp({\ve^2\, \pal_1 J_2})\, u^i
\eeq
leaves the form of the first Poisson structure unchanged while transforms the second one to the form
\beq
\exp({-\ve^2\, \ad_{\pal_1 J_2}})(Q_0+\sum_{k\ge 2} \ve^k \pal_1 {\tilde X}_k)
=Q_0+\sum_{k\ge 3} \ve^k \pal_1 {\bar X}_k
\eeq
Here the vector fields  $\bar X_{k}$ have the expressions
\beq
{\bar X}_k=\sum_{l=0}^{[k/2]} (-1)^l
 \frac{\left(\ad_{\pal_1 J_2}\right)^l}{l!} {\tilde X}_{k-2 l},\quad k\ge 3.
\eeq
and their components are homogeneous almost differential polynomials of degree $k$ which
depend only on finitely many of the jet coordinates. Here ${\tilde X}_1$ is assumed to be zero.
By repeatedly using the above procedure, we see that the deformation part
of (\ref{eq-2}) can be absorbed by a series of quasi-Miura transformations. So it is quasitrivial,
and the Theorem is proved.
\epf

{}From the proof of the Quasitriviality Theorem we see that the reducing transformation
\beq
u^i\mapsto u^i+\sum_{k\ge 1}\ve^k F_k(u,\dots,u^{(m_k)})
\eeq
of the bihamiltonian structure (\ref{eq-2}) has the properties that $F_k$ are homogeneous
almost differential polynomials of degree $k$ and that $m_k\le k$ when $k$ is odd and $m_k\le \frac32 k$
when $k$ is even.

\section{Reducing bihamiltonian PDEs}\label{sec-3}

In this section we study properties of the bihamiltonian systems (\ref{bih1}).
\begin{lem}\label{lem-f}
Let $I, J$ be two local functionals
\beq
I=\int p(u,u_x,\dots,u^{(N)})\, dx,\quad J=\int q(u,u_x,\dots,u^{(N)})\, dx
\eeq
that satisfy the relation
\beq
\pal_1 I=\pal_2 J.
\eeq
Then up to additions of total $x$-derivatives,
the densities $p, q$ do not depend on the jet coordinates $u_x,\dots, u^{(N)}$.
\end{lem}

\pf
Denote by $X$ the vector field $\pal_1 I-\pal_2 J$ and by $X^i$ its
components. Then from vanishing of
\beq
\frac{\pal X^i}{\pal u^{j,2 N+1}}=(-1)^{N}\left(f^i \frac{\pal^2 p}
{\pal u^{i,N}\pal u^{j,N}}-
g^i \frac{\pal^2 q}{\pal u^{i,N}\pal u^{j,N}}\right),
\eeq
we see that the functions $p, q$ satisfy
\beq
(u^i-u^j) \frac{\pal^2 q}{\pal u^{i,N}\pal u^{j,N}}=0,
\quad
\frac{\pal^2 p}{\pal u^{i,N}\pal u^{j,N}}-u^i \frac{\pal^2 q}{\pal u^{i,N}
\pal u^{j,N}}=0.
\eeq
So the functions $p, q$ can be represented by some functions $p_i, r_i, s$ as
\eqa
&&p=\sum_{i=1}^n p_i(u,\dots,u^{(N-1)},u^{i,N}),\nn\\
&&q=\sum_{i=1}^n \left[u^i p_i(u,\dots,u^{(N-1)},u^{i,N})+r_i(u,\dots,u^{(N-1)}) u^{i,N}+s(u,\dots,u^{(N-1)})\right]
\eeqa
By substituting these expressions of $p, q$ into the equations
\beq
\frac{\pal X^i}{\pal u^{i,2N}}=0
\eeq
we deduce that
\beq
f^i u^i \frac{\pal^2 p_i}{\pal u^{i,N}\pal u^{i,N}}=0
\eeq
Thus we can rewrite the functions  $p,q$ into the form
\eqa
p=\sum_{i=1}^n a_i(u,\dots,u^{(N-1)}) u^{i,N}+c(u,\dots,u^{(N-1)}),\nn\\
q=\sum_{i=1}^n b_i(u,\dots,u^{(N-1)}) u^{i,N}+d(u,\dots,u^{(N-1)}),
\eeqa
{}From the identity $\frac{\pal X^i}{\pal u^{j,2N}}=0$ for $i\ne j$ we have
\beq
\left(\frac{\pal^2 p}{\pal u^{i,N}\pal u^{j,N-1}}-\frac{\pal^2 p}{\pal u^{j,N}\pal u^{i,N-1}}\right)
=u^i \left(\frac{\pal^2 q}{\pal u^{i,N}\pal u^{j,N-1}}-\frac{\pal^2 q}{\pal u^{j,N}\pal u^{i,N-1}}\right)
\eeq
These equations imply that the functions $a_i, b_i$ satisfy
\beq
\frac{\pal a_i}{\pal u^{j,N-1}}=\frac{\pal a_j}{\pal u^{i,N-1}},\quad
\frac{\pal b_i}{\pal u^{j,N-1}}=\frac{\pal b_j}{\pal u^{i,N-1}}.
\eeq
So there exist functions $A(u,\dots,u^{(N-1)}), B(u,\dots,u^{(N-1)})$ such that
\beq
a_i=\frac{\pal A}{\pal u^{i,N-1}},\quad b_i=\frac{\pal B}{\pal u^{i,N-1}},\quad
i=1,\dots, n.
\eeq
Now we can replace the densities $p, q$ of the Hamiltonians $I, J$ respectively with
\beq
\tilde p= p-\pal_x A,\quad \tilde q=q-\pal_x B
\eeq
Then the new densities become independent of the jet variables $u^{i,N}, i=1,\dots,N$.
Repeating the above procedure successively, we arrive at the result of the lemma.
\epf

\noindent{\em Proof of Corollary \ref{thm2}.}
Let us assume that after the quasi-Miura transformation the Hamiltonians of the systems (\ref{bih1}) have the
expansion in $\ve$:
\beq
H_a=\sum_{k\ge 0} \ve^k H^{[k]}_a=\sum_{k\ge 0} \ve^k \int \tilde h^{[k]}_a(u,u_x,\dots,u^{(m_k)}) dx,
\quad a=1,2.
\eeq
Here $m_k$ are some positive integers which may also depend on the index $a$, and the functions
$h^{[k]}_a$ have degrees $k$. Due to the bihamiltonian property
\beq
\pal_1 H_1=\pal_2 H_2
\eeq
we know in particular that
\beq
\pal_1 H^{[1]}_1=\pal_2 H^{[1]}_2.
\eeq
Then the result of the above lemma implies that $H_1^{[1]}=H^{[1]}_2=0$. Similarly, we prove that all other
Hamiltonians $H^{[k]}_a,\ k\ge 2$ are trivial. The Theorem is proved. \epf

\begin{cor}
Any two bihamiltonian flows of the form (\ref{bih1}) that correspond to the same bihamiltonian structure
(\ref{eq-2}) mutually commute.
\end{cor}
\prf Denote by $X, Y$ the vector fields corresponding to the given bihamiltonian systems, then their commutator
$[X, Y]$ is also a bihamiltonian vector field of degree greater than $1$.
{}From Lemma \ref{lem-f} it follows that
under the quasi-Miura transformation reducing the bihamiltonian structure (\ref{eq-2}) to (\ref{eq-1})
this vector field must vanishes. Thus we proved the corollary.\epf

\section{Central invariants of bihamiltonian structures}\label{sec-4}

One of the important applications of the property of quasitriviality is the classification
of deformations of a given bihamiltonian structure of hydrodynamic type. The problem of classification
of quasitrivial infinitesimal deformations was solved in  \cite{LZ1}. It was also conjectured that
all deformations of the form (\ref{eq-2}) have reducing transformations. 
The Quasitriviality Theorem proves this
conjecture. In this section we reformulate  the main result of \cite{LZ1} in
order to describe the complete list of invariants of a bihamiltonian structure
with a given leading order $\{~,~\}^{[0]}_{1,2}$ modulo Miura-type
transformations (\ref{triv}). Recall that these transformations must depend
polynomially on the derivatives in every order in $\ve$.

Let us rewrite the bihamiltonian structure (\ref{eq-2}) in terms of the canonical coordinates
$u^i=u^i(w),\ i=1,\dots,n$
\begin{eqnarray}
&&\{u^i(x),u^j(y)\}_{a}=\{u^i(x),u^j(y)\}^{[0]}_a
\nn\\&&+\sum_{m\ge 1}\sum_{l=0}^{m+1} \ve^m A^{ij}_{m,l;a}(u;u_x,\dots,u^{(m+1-l)}) \delta^{(l)}(x-y),
\  a=1,2,\label{eq-2c}.
\end{eqnarray}
Then the functions $P^{ij}_a, Q^{ij}_a$ defined in (\ref{ss-31}) have the expressions
\beq\label{ss-31b}
P^{ij}_a(u)=A^{ij}_{1,2;a}(u),\ Q^{ij}_a(u)=A^{ij}_{2,3;a}(u),\ i,j=1,\dots,n,
a=1,2.
\eeq

\noindent{\em{Proof of Corollary \ref{thm51}}}.\
Let us first assume that the bihamiltonian structure (\ref{eq-2c}) has the following special form
\beq\label{ss-30}
(\varpi_1, \varpi_2+\ve^2 \gamma+{\mathcal O}(\ve^3)).
\eeq
Here $(\varpi_1,\varpi_2)$ denote the bihamiltonian structure given by the leading terms of (\ref{eq-2c}),
and $\gamma$ is a bivector which can be
represented as $\gamma=\pal_1 X$ through a vector
field with components that are homogenous differential polynomials of degree $2$.
Due to Theorem \ref{tmd-d2}, the vector field $X$ can be represented up to a Miura-type transformation in the
form
\beq
X=\pal_1 I-\pal_2 J
\eeq
where the functionals $I, J$ are defined by
\beq
I=\int \sum_{i=1}^n u^i \hat c_i(u^i) u^i_x \log u^i_x dx,\quad
J=\int \sum_{i=1}^n \hat c_i(u^i) u^i_x \log u^i_x dx
\eeq
with
\beq\label{ss-37}
\hat c_i(u^i)=\frac1{3 (f^i(u))^2} (\gamma)^{ii}_3,\quad i=1,\dots,n.
\eeq
Here $(\gamma)^{ii}_3$ denote the coefficients of $\delta'''(x-y)$ in the components $(\gamma)^{ii}$
of the bivector $\gamma$.
The main result of \cite{LZ1} together with the 
Quasitriviality Theorem
shows that any two bihamiltonian structures of the form (\ref{ss-30})
are equivalent {\em iff} they correspond
to the same set of functions $\hat c_i, i=1,\dots,n$. In the case that the two bihamiltonian structures of the
present theorem have the form (\ref{ss-30}), it is easy to see that $\hat c_i(u^i)=c_i(u)$, and the result
of the theorem follows.

Now return to the general form (\ref{eq-2c}) of the a bihamiltonian structure.
We redenote it as
\beq\label{ss-34}
(\varpi_1+\ve \al_1+\ve^2 \beta_1, \varpi_2+\ve \al_2+\ve^2 \beta_2)+{\mathcal O}(\ve^3).
\eeq
By using the result of Theorem \ref{lem-4} we can eliminate the linear in $\ve$ terms by a Miura-type transformation
\beq\label{ss-33}
u^i\mapsto \exp(-\ve X) u^i,\quad i=1,\dots, n
\eeq
given by a local vector field $X$ with components of the form
\beq
X^i=\sum_{j=1}^n X^i_j(u) u^j_x
\eeq
This implies that $P_1=\pal_1 X,\ P_2=\pal_2 X$ and this in turn yields
\beq
P_1^{ik}=-f^k(u) X^i_k(u)+f^i(u) X^k_i(u),\quad
P_2^{ik}=-g^k(u) X^i_k(u)+g^i(u) X^k_i(u)
\eeq
where the functions $P_1^{ij}, P_2^{ij}$ are defined by (\ref{ss-31b}). Solving the above system
we obtain
\beq\label{ss-32}
X^i_k(u)=\frac{P_2^{ki}-u^i P_1^{ki}}{f^k(u) (u^k-u^i)},\quad k\ne i.
\eeq

After the Miura-type transformation (\ref{ss-33}), the bihamiltonian structure (\ref{ss-34}) becomes
\beq\label{ss-36}
(\varpi_1+\ve^2 (\beta_1-\frac12 [X,\al_1]), \varpi_2+\ve^2 (\beta_2-\frac12 [X,\al_2]))
+{\mathcal O}(\ve^3).
\eeq
Then there exists a local vector field $Y$ such that
\beq
\pal_1 Y=\beta_1-\frac12 [X,\al_1].
\eeq
So the Miura-type transformation
\beq\label{ss-35}
u^i\mapsto \exp(-\ve Y) u^i,\quad i=1,\dots, n
\eeq
reduces the bihamiltonian structure (\ref{ss-36}) to the form of (\ref{ss-33})
\beq
(\varpi_1, \varpi_2+\ve^2 (\beta_2-\frac12 [X,\al_2]-\pal_2 Y))
+{\mathcal O}(\ve^3).
\eeq
and we need to compute the coefficients $\chi_3^{ij}$ of $\delta'''(x-y)$ in the components of the
bivector $\chi=\beta_2-\frac12 [X,\al_2]-\pal_2 Y$. By using the notations introduced in (\ref{ss-31b}) we have
\eqa
&&[X,\al_1]^{ii}_3=\sum_{k} X^i_k (P^{ik}_1-P^{ki}_1)=-2 \sum_{k\ne i} X^i_k P^{ki}_1,\nn\\
&&[X,\al_2]^{ii}_3=\sum_{k} X^i_k (P^{ik}_2-P^{ki}_2)=-2 \sum_{k\ne i} X^i_k P^{ki}_2,\nn\\
&&(\pal_2 Y)^{ii}_3=u^i (\pal_1 Y)^{ii}_3=u^i (\beta_1-\frac12 [X,\al_1])^{ii}_3=
u^i Q^{ii}_1-\frac12 u^i [X,\al_1]^{ii}_3
\eeqa
Here as above for any bivector $\eta$ we denote by $\eta^{ij}_3$ the
coefficient of $\delta'''(x-y)$ in the components $\eta^{ij}$.
These formulae together with the expressions (\ref{ss-32}) for $X^i_k, k\ne i$ show that the
functions $\hat c_i(u^i)$ that are defined by (\ref{ss-37}) with $\gamma$ replaced by $\chi$
coincide with the functions $c_i(u)$ introduced in (\ref{fc}). Thus we proved the theorem.
\epf

{}From this theorem it also easily follows the following corollary:
\begin{cor}
Any deformation (\ref{eq-2}) of the bihamiltonian structure (\ref{eq-1}) is equivalent,
under a appropriate Miura-type transformation, to a deformation  in which
only even powers of $\ve$ appear.
\end{cor}
This result can also be seen from the construction of the functionals $I, J$ in the proof of Theorem \ref{thm3}
and the argument given in the proof of the Quasitriviality Theorem.

\begin{thm}\label{thm53}
If we choose another representative
\beq\label{ss-40}
\{\, ,\, \}_1^{\tilde{}}=c  \{\, ,\, \}_2+d \{\, ,\, \}_1,\quad
\{\, ,\, \}_2^{\tilde{}}=a  \{\, ,\, \}_2+b \{\, ,\, \}_1,\quad ad-bc\ne 0
\eeq
of the bihamiltonian structure (\ref{eq-2c}), then the functions $c_i(u)$ that we define in (\ref{fc})
are changed to
\beq
\tilde c_i(\tilde u^i)=\frac{c u^i+d}{ad-bc} c_i(u^i),\quad i=1,\dots,n.
\eeq
where
\beq\label{ss-41}
\tilde u^i=\frac{a u^i+b}{c u^i+d},\quad i=1,\dots,n
\eeq
are the canonical coordinates of the bihamiltonian structure (\ref{eq-2c}) with respect to the new representative
(\ref{ss-40}).
\end{thm}

\begin{prf}
The result of the theorem is obtained by a straightforward calculation with the help of the formula
(\ref{ss-41}) and the tensor rule abided by $P_a^{ij}, Q_a^{ij}$ under the change of coordinates
$u^i\mapsto \tilde u^i(u)$.
\end{prf}

{}From the above theorem we see that, for the bihamiltonian structure
(\ref{eq-2c}), the following  $1/2$-forms
\beq
\Omega_i=c_i(u^i) (d u^i)^{\frac12}
\eeq
are invariant, up to a permutation, under the change of representative
(\ref{ss-40}) with $\begin{pmatrix} a&b\\ c& d\end{pmatrix}\in SL(2,\mathbb C)$.

\section{Examples and concluding remarks.}
Let us give some examples of bihamiltonian structures, 
their central invariants and reducing transformations.

\noindent{{\em Example 1}. {The Bihamiltonian structure of the KdV hierarchy \cite{gardner, magri, ZF}}} has the form
\begin{eqnarray}\label{kdv-pen}
&&\{w(x),w(y)\}_1=\delta'(x-y),\nn\\
&&\{w(x),w(y)\}_2=w(x)\delta'(x-y)+\frac12\,w'(x)\delta(x-y)+3\, c\,{\ve^2}\,\delta'''(x-y),\label{z2b}
\end{eqnarray}
The canonical coordinate is $u=w$, and the constant $c$ is the central
invariant. Up to terms of the order $\ve^6$ the reducing transformation
\cite{DZ1} is given\footnote{The
(inverse to) the reducing transformation for the KdV equation was constructed 
in \cite{bgi}. However, the action of this transformation onto the Poisson
pencil was not studied.}  
\begin{equation}
w=v +c\, {\epsilon^2} \p_x^2 \left( \log v_x\right)
+\frac{c^2\epsilon^4}{10} \p_x^2 \left( 5\frac{v^{(4)}}{ v_x^2}
- 21\frac{ v_{xx} v_{xxx}}{ v_x^3}
+16 \frac{v_{xx}^3}{ v_x^4}\right) + O(\epsilon^6).
\end{equation}
The Poisson pencils (\ref{kdv-pen}) with different values of $c$ are 
inequivalent with respect to Miura-type transformations.

\noindent{{\em {Example 2}}. {The Bihamiltonian structure related to the Camassa-Holm hierarchy \cite{CH, CHH, Fo, Fu, FF}}}
has the expression
\begin{eqnarray}
&&\{w(x),w(y)\}_1=\delta'(x-y)-\frac{\ve^2}8 \delta'''(x-y),\nn\\
&&\{w(x),w(y)\}_2=w(x)\delta'(x-y)+\frac12\,w'(x)\delta(x-y).\label{z9}
\end{eqnarray}
The canonical coordinate $u$ also coincides with $w$, and the central invariant $c(u)=\frac1{24} u$. As it was shown in
\cite{Lo} the reducing transformation is, up to $\ve^4$, given by
\eqa
&&w= v  + \epsilon^2\,
\pal_x \!\!\left(\frac{v \,v_{xx} }{24\,v_x }- \frac{v_x }{48} \right) \nn\\
&& \quad+
  \epsilon^4\pal_x\!\!\left(\frac{7\,v_{xx}^2}{2880\,v_x} +
  \frac{v\,v_{xx}^3}{180\,v_x^3} -
  \frac{{v}^2\,v_{xx}^4}{90\,v_x^5} - \frac{v_{xxx}}{512} -
  \frac{59\,v\,v_{xx}\,v_{xxx}}{5760\,v_x^2}\right.\nn\\
  &&\qquad\left. +
  \frac{37\,{v}^2\,v_{xx}^2\,v_{xxx}}{1920\,v_x^4} -
  \frac{7\,{v}^2\,v_{xxx}^2}{1920\,v_x^3} +
  \frac{5\,v\,v^{(4)}}{1152\,v_x} -
  \frac{31\,{v}^2\,v_{xx}\,v^{(4)}}{5760\,v_x^3} +
  \frac{{v}^2\,v^{(5)}}{1152\,v_x^2}\right).
\eeqa



\noindent{{\em {Example 3}}. The Bihamiltonian structure related to the multi-component
KdV-CH (Camassa-Holm) hierarchy.
Define
\beq
\mathcal{D}_i=w^i \delta'(x-y)+\frac{w^i_x}2 \delta(x-y)+ a_i\frac{\ve^2}8
\delta'''(x-y), \quad i=0,1,\dots,n
\eeq
Here $w^0=1$ and $a_0,\dots,a_n$ are given constants with at least one
nonzero. Define also the numbers
\beq
f_m^{ij}=\begin{cases}-1& i, j\le m\\ +1 & i,j\ge m+1\\ 0 &\text{otherwise}\end{cases}
\eeq
Then we have the following $n+1$ compatible Hamiltonian structures
\beq\label{ss-50}
\{w^i(x), w^j(y)\}_m=(-1)^m f^{ij}_m\, {\mathcal D}_{i+j-m-1},\quad 1\le i, j\le n,\ m=0,1,\dots,n.
\eeq
When $i<0$ or $i>n$ we assume ${\mathcal D}_i=0$.
These Hamiltonian structures were introduced\footnote{To our best knowledge,
connections of these bihamiltonian structures with the Camassa - Holm equation
and its multicomponent generalizations was never considered in the literature.} in the study of the hierarchies
of integrable systems (called the coupled KdV hierarchies) associated with the compatibility conditions of the linear systems of the form \cite{fordy, fordy2, fordy3, alonso}
\eqa
&&\left(\frac12 \left(\ve\,\pal_x\right)^2+A(w;\lm)\right) \psi=0,\\
&&\psi_t=\frac12B \psi_x-\frac14 B_x \psi.
\eeqa
Here $A(w;\lm)$ has the expression
\beq
A(w;\lm)=\frac{\sum_{i=0}^n(-1)^iw^i\lm^{n-i}}{\sum_{i=0}^n(-1)^ia_i\lm^{n-i}}
\eeq
and $B$ is certain polynomial or Laurent polynomial in $\lm$ 
with coefficients that are
differential polynomials of $w^1,\dots, w^n$ 
which can be chosen according to the equation
\beq
A_t=A B_x+\frac{A_x}2B+\frac{\ve^2}8B_{xxx}.
\eeq
As it was shown by Ferapontov in \cite{fev2}, if a system of hydrodynamic type with $n$ depend variables
possesses $n+1$ compatible Hamiltonian structures of hydrodynamic type, then this $n+1$-Hamiltonian
structure must be equivalent to the one obtained from the leading terms of (\ref{ss-50}).

{}From (\ref{ss-50}) we readily have the following bihamiltonian
structures
\beq
{\mathcal B}_{k,l}=(\{\,,\,\}_k, \{\,,\,\}_l),\quad k \ne l.
\eeq
Denote by $\lm_1(w),\dots, \lm_n(w)$ the roots of the polynomial $P(\lm)=\lm^n+w^1 \lm^{n-1}+\dots+w^n$.
Then the canonical coordinates for the bihamiltonian structure ${\mathcal B}_{k,l}$
are given by $u^i=(\lm_i)^{k-l},\ i=1,\dots,n$, and the central invariants $c_1,\dots, c_n$ have the expressions
\beq\label{nCHdfm}
c_i(u^i)=\frac{\sum_{i=0}^n(-1)^ia_i(\lm_i)^{n-i}}{24\, (l-k) (\lm_i)^{n-1-l}},\quad i=1,\dots,n.
\eeq

In particular, for the one-component case $n=1$, choosing $a_0=0, a_1=1$ we get the bihamiltonian structure
${\mathcal B}_{1,0}$ which coincides with (\ref{z2b}) for the KdV hierarchy.
The choice $a_0=-1, a_1=0$ yields the bihamiltonian
structure (\ref{z9}) of the Camassa-Holm hierarchy. In general, we call
the hierarchy generated by the bihamiltonian
structure ${\mathcal B}_{k,l}$ the {\it multi-component KdV-CH hierarchy}.


For the case when $n=2$, the above defined bihamiltonian structure ${\mathcal B}_{2,1}$ yields,
with different choices of the constants $a_0,a_1,a_2$ and up to certain Miura-type transformation and rescaling of
the Poisson structures, the following four bihamiltonian structures that
appeared in the literature. They have the same leading terms
\eqa
&&\{\vp(x),\vp(y)\}_1^{[0]}=0,\ \{\rho(x),\vp(y)\}_1^{[0]}=\dl'(x-y),\ \{\rho(x),\rho(y)\}_1^{[0]}=0, \nn \\
&&\{\vp(x),\vp(y)\}_2^{[0]}=2\,\dl'(x-y),\ \{\rho(x),\vp(y)\}_2^{[0]}=\vp(x)\,\dl'(x-y), \nn \\
&&\{\rho(x),\rho(y)\}_2^{[0]}=2\rho(x)\,\dl'(x-y)+\rho'(x)\,\dl(x-y).\label{lt}
\eeqa
A bihamiltonian structure related to the nonlinear Schr\"odinger hierarchy is given by the above
brackets with the only difference \cite{bonora, LZ1}
\beq
\{\rho(x),\vp(y)\}_2=\{\rho(x),\vp(y)\}_2^{[0]}+\ve \dl''(x-y).
\eeq
After the Miura-type transformation
\beq
w^1=2\vp,\quad w^2=\vp^2-4\rho+2 \ve \vp_x \\
\eeq
it is transformed
to the bihamiltonian structure $8\, {\mathcal B}_{2,1}$ with the choice of constants
$a_0=a_1=0, a_2=-8$. Here $8\, {\mathcal B}_{2,1}$ denotes the bihamiltonian structure obtained from
${\mathcal B}_{2,1}$ by the mutiplication of a overall factor $8$.


In \cite{LZ1} a generalization of the Camassa-Holm hierarchy is introduced which is called
the 2-component Camassa-Holm hierarchy. It is reduced to the usual Camassa-Holm hierarchy
under a natural constraint on its two dependent variables. The related bihamiltonian structure
is defined by the brackets (\ref{lt}) except
\beq
\{\rho(x),\vp(y)\}_1=\{\rho(x),\vp(y)\}_1^{[0]}+\ve \dl''(x-y).
\eeq
After the Miura-type transformation
\beq
w^1=2\vp+2\ve\vp_x,\quad w^2=\vp^2-4\rho
\eeq
it is converted, up to the approximation to $\ve^2$, to the bihamiltonian structure $8\, {\mathcal B}_{2,1}$ with the choice of constants
$a_0=-8, a_1=0, a_2=0$.

In \cite{kersten} the bihamiltonian structure for the so called
classical Boussinesq hierarchy is given.
It is defined by the brackets (\ref{lt}) except for
\beq
\{\rho(x),\rho(y)\}_2=\{\rho(x),\rho(y)\}_2^{[0]}+\frac12 \ve^2 \dl'''(x-y).
\eeq
After the Miura-type transformation
\beq
w^1=2\vp,\quad w^2=\vp^2-4\rho
\eeq
it is transformed to the bihamiltonian structure $8\, {\mathcal B}_{2,1}$ with the choice of constants
$a_0=0, a_1=0, a_2=-8$.

Note that, for the bihamiltonian structure related to the nonlinear
Schr\"odinger hierarchy,
 by moving the perturbation term from the second Poisson bracket to
the first one we
obtain the bihamiltonian structure of the 2-component Camassa-Holm hierarchy.
Doing precisely in the same way we obtain from the
above bihamiltonian structure of the classical Boussinesq hierarchy the one
that is defined by (\ref{lt}) except  for the bracket
\beq\label{fpb}
\{\rho(x),\rho(y)\}_1=\{\rho(x),\rho(y)\}_1^{[0]}-\frac12 \ve^2 \dl'''(x-y).
\eeq
After the change of dependent variables
\beq
w^1=2\vp,\quad w^2=\vp^2-4\rho
\eeq
it is transformed to the bihamiltonian structure $8\, {\mathcal B}_{2,1}$ with the choice of constants
$a_0=0, a_1=8, a_2=0$. This bihamiltonian structure is related to the Ito type equations \cite{ito, kup}.

The bihamiltonian structures related to the nonlinear Schr\"odinger hierarchy and
the classical Boussinesq hierarchy are equivalent. Indeed, their central invariants are given by $c_1=c_2=\frac1{24}$.
The central invariants for the bihamiltonian structure related to the 2-component Camassa-Holm
hierarchy are given by $c_1=\frac{(u^1)^2}{24},\ c_2=\frac{(u^2)^2}{24}$,
and those for the bihamiltonian
structure defined by (\ref{lt}) and (\ref{fpb}) have the form $c_1=\frac{u^1}{24},\ c_2=\frac{u^2}{24}$.

We omit here the presentation of the reducing transformations of the above bihamiltonian structures due to
their cumbersome expressions.

\noindent{\em Example 4}. The 
equations of motion of one-dimensional isentropic gas with the equation of
state $p=\frac{\kappa}{\kappa+1}\rho^{\kappa+1}$ read 
\beq\label{fluid}
u_t+\left( \frac{u^2}{2} + \rho^\kappa\right)_x=0,
\quad
\rho_t +(\rho\,u)_x =0.
\eeq
Here $\kappa$ is an arbitrary parameter, $\kappa\neq 0,
\, -1$. For a gas with
$m$ degrees of freedom one has
$$
\kappa = \frac2{m}
$$
(see, e.g., \cite{daf}). This is a weakly simmetrizable system with
$$
\eta=\left(\begin{array}{cc} 0 & 1 \\ 1 & 0 \end{array}\right).
$$
This gives the first Poisson structure of the equations with the Poisson
brackets
\beq\label{I2-pb0}
\left\{ u(x), \rho(y)\right\}^{[0]}_1 =\delta'(x-y),
\eeq
other brackets vanish.
The second hamiltonian structure  
\eqa\label{I2-pb}
&&\{u(x),u(y)\}^{[0]}_2=2\rho^{\kappa-1}(x)\,\delta'(x-y)
+ \rho^{\kappa-1}_x\,\delta(x-y),\nn\\
&&\{u(x),\rho(y)\}^{[0]}_2=u(x)\,\delta'(x-y)+\frac1{\kappa}\,u'(x)\,
\delta(x-y),\nn\\
&&\{\rho(x),\rho(y)\}^{[0]}_2=\frac1{\kappa}\left(2\,\rho(x)\,\delta'(x-y)+
\rho'(x)\,\delta(x-y)\right).
\eeqa
was found in \cite{nutku, nutku-olver}. 

As it was shown in \cite{DZ1} the isentropic gas equations
have the following deformation which preserves the bihamiltonian property
(up to corrections\footnote{In principle one can continue the expansions till an
arbitrary order in $\ve$. However, the computations become very involved.} of order
$\ve^6$):
\eqa\label{poly-def}
&&\frac{\partial u}{\partial t}+\partial_x\left\{
\frac{u^2}{2}+{\rho^\kappa}
+\epsilon^2\left[
\frac{\kappa(\kappa-2)}{8}\rho^{\kappa-3}\,\rho_x^2
+\frac{\kappa^2}{12}\,\rho^{\kappa-2}\rho_{xx}\right]\right.
\nn\\
&&\quad +\epsilon^4 (\kappa-2)(\kappa-3)\left[
a_{1}\, \rho^{-4} \,u_x^2\,\rho_x^2+
        a_{2}\, \rho^{\kappa-6} \,\rho_x^4+
        a_{3}\, \rho^{-3} \,u_{xx}\,u_x\,\rho_x\right.
\nn\\
&&\qquad+
        a_{4}\, \rho^{-2} \,u_{xx}^2+
        a_{5}\, \rho^{-3} \,u_x^2\,\rho_{xx}
+
        a_{6}\, \rho^{\kappa-5} \,\rho_x^2\,\rho_{xx}+
        a_{7}\, \rho^{\kappa-4} \,\rho_{xx}^2+
        a_{8}\, \rho^{-2} \,u_x\,u_{xxx}
\nn\\
&&\qquad\quad\left. \left.+
        a_{9}\, \rho^{\kappa-4} \,\rho_x\,\rho_{xxx}+a_{10} \rho^{-\kappa-2} u_x^4\right]+
        \epsilon^4 \frac{\kappa\,(\kappa^2-4)}{360} \rho^{\kappa-3}
        \,\rho_{xxxx}\right\} ={\mathcal O}(\epsilon^6),
\nn
\eeqa
\eqa
&&\frac{\partial \rho}{\partial t}+\partial_x\left\{
\rho\, u+\epsilon^2\,\left(\frac{(2-\kappa)(\kappa-3)}{12\,\kappa\,\rho}\,u_x\,\rho_x+
\frac16\,u_{xx}\right)\right.
\nn\\
&&\nn\\&&
\quad+\epsilon^4 (\kappa-2)(\kappa-3)\left[
b_{1}\, \rho^{-4} \,u_x\,\rho_x^3+
        b_{2}\, \rho^{-3} \,\rho_x^2\,u_{xx}+
        b_{3}\, \rho^{-3} \,u_x\,\rho_x\,\rho_{xx}\right.
\nn\\
&&\qquad+
        b_{4}\, \rho^{-2} \,u_{xx}\,\rho_{xx}+
        b_{5}\, \rho^{-2} \,u_{xxx}\,\rho_x+
        b_{6}\, \rho^{-2} \,u_x\,\rho_{xxx}+
        b_{7}\, \rho^{-1} \,u_{xxxx}\nn\\
&&\left.\left.\qquad+ b_8\, \rho^{-\kappa-1} u_x^2 u_{xx}+b_9\, \rho^{-\kappa-2} u_x^3 \rho_x\right]\right\}
={\mathcal O}(\epsilon^6)
\eeqa
The coefficients are given by
\eqa
&&a_{1}={\frac{
     18 + 75\,\kappa - 15\,{\kappa^2} + 20\,{\kappa^3}+2\kappa^4  }{2880\,{\kappa^3}}},
~~a_{2}={\frac{
     6 + 113\,\kappa + 409\,{\kappa^2} - 185\,{\kappa^3} + 17\,{\kappa^4}  }{5760\,
     {\kappa^2}}}
\nn\\&&\nn\\
&&a_{3}=-{\frac{
     18 +11\,\kappa +3\,{\kappa^2}}{720\,{\kappa^2}}},~~
a_{4}={\frac{7}{720\,{\kappa}}},~~
a_{5}={\frac{
      -6 + 3\,\kappa - {\kappa^2}  }{480\,{\kappa^2}}}
\nn\\&&\nn\\
&&a_{6}={\frac{
      -6 - 39\,\kappa - 10\,{\kappa^2} +5\,{\kappa^3}}{480\,
     {\kappa}}},~~
a_{7}=\frac{
      14 + 5\,\kappa + 5\,{\kappa^2} }{1440},
\quad
a_{8}={\frac{1 }{120\,\kappa}}\nn\\
&&a_{9}={\frac{ 2 + 5\,\kappa  }
     {240}},\quad
a_{10}=-\frac{(\kappa+2)(\kappa+3)(\kappa^2-1)}{5760\,\kappa^4}.
     \eeqa

\eqa
&&b_{1}={\frac{
      42 + 83\,\kappa - 53\,{\kappa^2} + 8\,{\kappa^3}
      }{1440\,{\kappa^3}}},~~
b_{2}=-{\frac{
     6 + 35\,\kappa -24\,{\kappa^2} +5\,{\kappa^3}  }{720\,{\kappa^3}}}
\nn\\&&\nn\\
&&b_{3}=-{\frac{
      12 + 40\,\kappa - 13\,{\kappa^2} + 5\,{\kappa^3}  }{720\,{\kappa^3}}},~
b_{4}={\frac{
      6 - 4\,\kappa + {\kappa^2}  }{180\,{\kappa^2}}},~
b_{5}={\frac{
      6 + \kappa + {\kappa^2} }{720\,{\kappa^2}}}\nn\\
&&\nn\\
&&b_{6}={\frac{
      6 + \kappa + {\kappa^2}  }{720\,{\kappa^2}}},~
b_{7}=-{\frac{1 }{360\,\kappa}},\quad
b_8=-\frac{(\kappa+2)(\kappa+3)}{720\,\kappa^4},\nn\\
&&b_9=\frac{(\kappa+1)(\kappa+2)(\kappa+3)}{1440\,\kappa^4}.
\eeqa
The corresponding bihamiltonian structure (at the approximation up to $\ve^4$)
is given in Section 4.2.3 of \cite{DZ1}, the central invariants
are $c_1=c_2=\frac1{24}$. The above system can be represented as
\eqa
&&\frac{\pal u}{\pal t}=\{H_1,u(x)\}_1=\frac{\kappa}{\kappa+1} \{H_2,u(x)\}_2,\nn\\
&&\frac{\pal \rho}{\pal t}=\{H_1,\rho(x)\}_1=\frac{\kappa}{\kappa+1}
\{H_2,\rho(x)\}_2.
\eeqa
Here the densities $h_1, h_2$ of the Hamiltonians $H_1, H_2$ have 
the expressions
\eqa
&&
h_1 =\frac12 \rho\, u^2 +{\rho^{\kappa +1}\over \kappa+1} +\Delta h_1
\nn\\
&&
h_2 = \rho\, u +\Delta h_2
\nn
\eeqa
where
\eqa
&&\Delta h_1=u\,\Delta h_2
-
\frac{\ve^2}{24\kappa}\,\left(({\kappa}^{2}-3\,\kappa+6)\,u_x^2+
\kappa(2\,{{\kappa}}^{2}-5\,\kappa+6)\,\rho^{\kappa-2}\,\rho_x^2\right)
\nn\\&&\quad +
\ve^4\,\frac{(\kappa-2)\,(\kappa-3)}{240\,\kappa^3}\left(
-\frac13\,\kappa\,\left({\kappa}^{2}-4\,\kappa+6\right)\,\rho^{-1}\,u_x\,u_{xxx}\right.\nn\\
&&\quad
+\frac13\,\kappa\,\left (2\,{\kappa}^{2}-13\,\kappa+12\right )\,\rho^{-2}\,u_x\,u_{xx}\,\rho_x
+{\frac {1}{72}}\,{\frac {\left (3\,\kappa+5\right )\left (\kappa+3
\right )\left (\kappa+2\right )}{\kappa}}\,\rho^{-\kappa-1}\,u_x^4 \nn \\
&&\quad-\frac1{12}\,
\left (2\,\kappa-3\right )\left (\kappa+3\right )\left ({
 \kappa}+2\right )\,\rho^{-3}\,u_x^2\,\rho_x^2+
{\frac {{\kappa}^{2}\left (\kappa-1\right )\left (3\,{{ \kappa
}}^{2}-8\,\kappa+12\right )}{2(\kappa-3)}}\,\rho^{\kappa-3}\,\rho_{xx}^2
\nn\\&&\quad\left.-{\frac1{72}}\,\kappa\,\left (\kappa-1\right )\left (16\,{{
\kappa}}^{4}-100\,{\kappa}^{3}+229\,{\kappa}^{2}-211\,\kappa+6
\right )\,\rho^{\kappa-5}\,\rho_x^4\right),
\eeqa
\eqa
&&\Delta h_{2}=-\ve^2\frac{(\kappa-2)\,(\kappa-3)}{12\,\kappa}\rho^{-1}\,u_x\,\rho_x
\nn\\&&\qquad
+\ve^4\frac{(\kappa-2)\,(\kappa-3)}{720\,\kappa^3}
\left[-2\,\kappa\,\left ({\kappa}^{2}-8\,\kappa+6\right )
\,\rho^{-2}\,u_x\,\rho_{xxx}\right.\nn\\&&\qquad
+\kappa\,\left (7\,{\kappa}^{2}-61\,\kappa+42\right )
\,\rho^{-3}\,u_x\,\rho_x\,\rho_{xx}+
(-5\,{\kappa}^{3}+{\frac {79}{2}}\,{\kappa}^{2}-{\frac {55}{2}}\,
\kappa+3)
\,\rho^{-4}\,u_x\,\rho_x^3\nn\\
&&\qquad\qquad+
\frac1{6\kappa}\,\left.\left (\kappa+3\right )\left (\kappa+2\right )(\kappa+1)
\,\rho^{-k-2}\,u_x^3\,\rho_x
\right].
\eeqa
To write down the reducing transformation of the perturbed system of the one dimensional
isentropic gas and its bihamiltonian structure, we introduce
the operators ${\mathcal T_1}, {\mathcal T_2}$
\eqa
&&{\mathcal T}_1\, u^{(m)}=u^{(m+1)},\quad {\mathcal T}_1\, \rho^{(m)}=\rho^{(m+1)},\nn\\
&&
{\mathcal T}_2\, u^{(m)}=\pal_x^m \left(\kappa \rho^{\kappa-2}\,\rho_x\right),\quad
{\mathcal T}_2\, \rho^{(m)}=u^{(m+1)},\ m\ge 0.
\eeqa
We will use Greek subscripts for the result of acting of the operators on the
functions $\rho$ and $u$, i.e.
$$
\rho_{\alpha_1 \alpha_2 \dots}:= 
{\mathcal T}_{\alpha_1} {\mathcal T}_{\alpha_2} \dots \rho, \quad 
u_{\alpha_1 \alpha_2 \dots}:= 
{\mathcal T}_{\alpha_1} {\mathcal T}_{\alpha_2} \dots u.
$$
Define the functions
\beq
{\mathcal F_1}=\frac1{24} \log\left(\kappa \rho^{\kappa-2} \rho_x^2-u_x^2\right)-\frac1{24}
\frac{(\kappa-2)(\kappa-3)}{\kappa} \log\rho,
\eeq
\eqa
&&{\cal F}_2=\frac1{1152}\,\rho_{\al_1\al_2\al_3\al_4}
M^{\al_1\al_2}\,
M^{\al_3\al_4}-\frac1{360}\,\rho_{\al_1\al_2\al_3}
\rho_{\al_4\al_5\al_6}
M^{\al_1\al_4}\,
M^{\al_2\al_5}\, M^{\al_3\al_6}\nn\\
&&-\frac1{1152}\,\rho_{\al_1\al_2}
\rho_{\al_3\al_4\al_5\al_6}M^{\al_1 \al_3}\, M^{\al_2\al_4}\, M^{\al_5\al_6}
\nn\\
&&+
\frac1{360}\,\rho_{\al_1\al_2}\rho_{\al_3\al_4\al_5}
\rho_{\al_6\al_7\al_8}
M^{\al_1\al_3}\,
M^{\al_2\al_6}\,M^{\al_4\al_7}\, M^{\al_5\al_8}\nn
\\&&+{(\kappa-2)(\kappa-3)}\,D^{-2}\,\left[-\frac1{240} \kappa\,\rho^{2\,k-5}\,\rho_{xxx}\,
\rho_x^3
+\frac{11}{2880} \kappa\,\rho^{2\,\kappa-5}\,\rho_{xx}^2\,\rho_x^2 \right.\nn\\
&&+
\left(-\frac7{5760} \kappa^2+\frac{19}{5760}\kappa+\frac7{960}\right)\,\rho^{2\,\kappa-6}\,\rho_{xx}\,\rho_x^4
+
\frac{11}{2880}\,\rho^{\kappa-3}\,\rho_x^2\,u_{xx}^2
\nn
\\&&
-\frac1{5760 \kappa}(\kappa^4-9\kappa^3+\kappa^2+53 \kappa+6)\,\rho^{2\,\kappa-7}\,\rho_x^6
+\frac1{240}\,\rho^{\kappa-3}\,\rho_x^2\,u_{xxx}\,u_x\nn\\
&& +
\frac1{240}\,\rho^{\kappa-3}\,\rho_x\,\rho_{xxx}\,u_x^2-\frac{11}{720}\,\rho^{\kappa-3}\,\rho_x\,u_x\,u_{xx}\,\rho_{xx} +
\frac{11}{2880}\,\rho^{\kappa-3}\,u_x^2\,\rho_{xx}^2\nn
\\&&-\frac1{1440} (11\kappa-21)\,\rho^{\kappa-4}\,\rho_x^3\,u_{xx}\,u_x +
\frac1{2880\kappa} (22\kappa^2-47\kappa-42)\,\rho^{\kappa-4}\,\rho_x^2\,u_x^2\,\rho_{xx}\nn
\\&&+\frac1{5760 \kappa^2} (12\kappa^4-45\kappa^3+15\kappa^2+101\kappa+6)\,\rho^{\kappa-5}\,u_x^2\,\rho_x^4 -
\frac{1}{240\kappa}\,\rho^{-1}\,u_{xxx}\,u_x^3\nn
\\&&+\frac{11}{2880\kappa}\,\rho^{-1}\,u_x^2\,u_{xx}^2 +
\frac1{1440\kappa}\,\rho^{-2}\,u_{xx}\,u_x^3\,\rho_x
+\frac1{5760\kappa^2} (7\kappa^2-13\kappa+42)\,\rho^{-2}\,\rho_{xx}\,u_x^4\nn\\
&& \left.-
\frac1{5760\kappa^3}(8\kappa^3-31\kappa^2+43\kappa-6)\,\rho^{-3}\,u_x^4\,\rho_x^2
-\frac1{5760\kappa^4} (\kappa+3)(\kappa+2)\,\rho^{-\kappa-1}\,u_{x}^6\,\right]
\nn
\eeqa
Here the matrix $M(M^{\al\beta})$ and differential polynomial $D$ read
\beq
M=D^{-1}\begin{pmatrix} -\kappa \rho^{\kappa-2}\,\rho_x & u_x  
\\ u_x &  -\rho_x \end{pmatrix},
\quad D= u_x^2 -\kappa\, \rho^{\kappa -2} \rho_x^2,
\eeq
Then the reducing transformation is given by the formula
\beq
u\mapsto u+{\mathcal T}_1{\mathcal T}_2\left(\ve^2 {\mathcal F}_1+ \ve^4 {\mathcal F}_2\right),\quad
\rho\mapsto \rho+{\mathcal T}_1{\mathcal T}_1\left(\ve^2 {\mathcal F}_1+ \ve^4 {\mathcal F}_2\right).
\eeq
We leave as an exercise for the reader to check that the denominator 
$D\neq 0$ on the monotone solutions.

\bigskip

In conclusion let us formulate some open problems.

{\bf Problem 1}. Study convergence of the reducing transformations for the case
of analytic in $\ve$ perturbations (\ref{pert4}) (e.g., for the case of polynomial
dependence on $\ve$).

\smallskip

{\bf Problem 2}. Are there more wide classes of perturbations of systems of
hyperbolic PDEs admitting reducing transformation? The natural candidate to be
considered is the perturbations of the so-called {\it semi-hamiltonian} systems
in the Tsarev's sense \cite{tsarev}, i.e., hyperbolic systems written in the
diagonal form and possessing a complete family of commuting flows.

\smallskip

{\bf Problem 3}. According to our results, classes of equivalence of semisimple
bihamiltonian structures depend at most on $n(n+1)$ arbitrary functions of one
variable. Prove existence of such bihamiltonian structures for an arbitrary
choice of these functional parameters.

\def\thetheorem{A.\arabic{theorem}}
\def\theprop{A.\arabic{prop}}
\def\thelemma{A.\arabic{lemma}}
\def\thecor{A.\arabic{cor}}
\def\theexam{A.\arabic{exam}}
\def\theremark{A.\arabic{remark}}
\def\theequation{A.\arabic{equation}}

\setcounter{equation}{0}
\appendix
\makeatletter
\renewcommand{\@seccntformat}[1]{{Appendix:}\hspace{-2.3cm}}
\makeatother
\renewcommand{\thesection}{Appendix:}
\section{\quad\qquad \ \ Bihamiltonian structures of hydrodynamic type}

In this Appendix we will describe in more details, following \cite{FEV}, the defining equations for semisimple bihamiltonian structures of hydrodynamic type as well as their Lax pair representation.

We will work in the canonical coordinates $u^1$, \dots, $u^n$ (see Lemma \ref{lem-fs} above). 
Introduce the classical {\it Lam\'e coefficients}
$$
H_i(u):= f_i^{-1/2}(u), \quad i=1, \dots, n
$$
and the {\it rotation coefficients}
\beq\label{rot}
\gamma_{ij}(u):= H_i^{-1} \pal_i H_j, \quad i\neq j.
\eeq
Here, as usual
$$
\pal_i =\frac{\pal}{\pal u^i},
$$
no summation over repeated indices will be assumed within this section.
The classical {\it Lam\'e equations}
\eqa\label{lame1}
&&
\pal_k \gamma_{ij} =\gamma_{ik} \gamma_{kj}, \quad i, \, j, \, k ~ \mbox{distinct}
\\
&&
\pal_i \gamma_{ij} +\pal_j\gamma_{ji} 
+\sum_{k\neq i, \, j} \gamma_{ki}\gamma_{kj}=0, \quad i\neq j
\label{lame2}
\eeqa
describe diagonal metrics of curvature zero\footnote{Integrability of the system
(\ref{lame1}) - (\ref{lame2}) was discovered by V.Zakharov\cite{zakh}.}. Adding the equations
\beq\label{lame3}
u^i \pal_i \gamma_{ij}+u^j\pal_j \gamma_{ji} 
+\sum_{k\neq i, \, j} u^k \gamma_{ki}\gamma_{kj} 
+\frac12 (\gamma_{ij}+\gamma_{ji})=0, \quad i\neq j
\eeq
one obtains the defining relations for semisimple Poisson pencils of hydrodynamic type.
The solutions to the system (\ref{lame1}) - (\ref{lame3}) are parametrized by $n(n-1)$
arbitrary functions of one variable. Indeed, one can freely choose the functions
$$
\gamma_{ij}(u^1_0, \dots, u^j, \dots, u^n_0)
$$
near a given point
\beq\label{po}
u_0=(u^1_0, \dots, u^n_0), \quad u^i_0\neq u^j_0, \quad u^i_0\neq 0.
\eeq

The equations (\ref{lame1}) - (\ref{lame3}) can be represented as the compatibility conditions of the linear system
\eqa\label{lax1}
&&
\pal_i\psi_j =\gamma_{ji} \psi_i, \quad i\neq j
\nn\\
&&
\pal_i \psi_i + \sum_{k\neq i} \gamma_{ki} \frac{u^k -\lambda}{u^i-\lambda} \psi_k
+\frac{1}{2\, (u^i-\lambda)} \psi_i=0.
\eeqa
(``Lax pair'' with the spectral parameter $\lambda$ for (\ref{lame1}) - (\ref{lame3})).
The solutions to the linear system (\ref{lax1}) are closely related with the common first integrals of the bihamiltonian systems of hydrodynamic type, i.e. with the Casimirs
of the Poisson pencil
\eqa\label{casi}
&&
\{~, I\}^{[0]}_2 - \lambda\{~,I\}^{[0]}_1 =0, \quad  I=\int P(u)\, dx
\nn\\
&&
\pal_i P(u) =\psi_i H_i, \quad i=1, \dots, n.
\eeqa
As we already know (see Lemma \ref{lem-fs1} above) the bihamiltonian systems
are all diagonal in the canonical coordinates
\beq\label{hs1}
u^i_t +V^i(u) u^i_x =0, \quad i=1, \dots, n.
\eeq
The characteristic velocities $V^i(u)$ are determined from the following linear system
\eqa\label{vel}
&&
\pal_k\chi_i =\gamma_{ki} \chi_k, \quad i\neq k
\\
&&
\chi_i = H_i V^i, \quad i=1, \dots, n.
\eeqa
For the given rotation coefficients $\gamma_{ij}(u)$ satisfying (\ref{lame1}) - (\ref{lame3}) the solutions to (\ref{vel}) can be reconstructed with ambiguity of $n$ arbitrary functions of one variable.
In particular, the Lame coefficients $\chi_i=H_i(u)$ give a solution to (\ref{vel}). They correspond to the spatial translations $V^i(u) \equiv 1$, $i=1, \dots, n$.

Finally, to reconstruct the flat pencil of metrics starting from given solution to (\ref{lame1}) - (\ref{lame3}) near a given point  (\ref{po}) one has to choose a solution $\chi_1(u)$, \dots, $\chi_n(u)$ such that
$$
\chi_i(u_0)\neq 0, \quad i=1, \dots, n.
$$
Then we put
\beq\label{pen1}
g^{ij}_1(u) = \chi_i^{-2}(u)\delta_{ij}, \quad g^{ij}_2(u) = u^i\chi_i^{-2}(u)\delta_{ij}.
\eeq
The flat coordinates of the metrics correspond to particular solutions of the system
(\ref{lax1}). Namely, to find flat coordinates for the first metric one has to choose a fundamental system of solutions
$$
\psi_i^\alpha(u), \quad \alpha =1, \dots, n, \quad \det \left( \psi_i^\alpha(u_0)\right) \neq 0
$$
to the following linear overdetermined system
\eqa\label{fl1}
&&
\pal_i \psi_j =\gamma_{ji}\psi_i, \quad i\neq j
\nn\\
&&
\pal_i\psi_i +\sum_{k\neq i} \gamma_{ki} \psi_k =\frac12 \psi_i
\eeqa
obtained from (\ref{lax1}) at $\lambda=\infty$. Then the flat coordinates $v^\alpha$
are defined by quadratures
\beq\label{fl2}
d v^\al =\sum_{i=1}^n  \chi_i \psi_i^\alpha du^i, \quad \alpha = 1, \dots, n.
\eeq
Flat coordinates for the second metric are constructed in a similar way by using a fundamental system of solutions to (\ref{lax1}) at $\lambda=0$.

We deduce that semisimple bihamiltonian structures of hydrodynamic type with $n$ dependent variables
are parametrized by $n^2$ arbitrary functions of one variable. For $n\leq 2$ the equations (\ref{lame1}) - (\ref{lame3}) are linear. So an
explicit parametrization of the Poisson pencils is available \cite{mokhov}. The equations become nonlinear starting from $n\geq 3$. All known so far nontrivial solutions
are obtained within the theory of Frobenius manifolds. In this case the rotations coefficients are symmetric
$$
\gamma_{ji}=\gamma_{ij}
$$
(the so-called {\it Egoroff metrics})
and the equations (\ref{lame1}) - (\ref{lame3}) are reduced to isomonodromy deformations \cite{Du1}. We will study more general case in subsequent publications.

\noindent{Emails:}\newline
 dubrovin@fm.sissa.it,\ lsq99@mails.tsinghua.edu.cn,\ youjin@mail.tsinghua.edu.cn}
\end{document}